\documentclass[10pt,a4paper]{article}
\usepackage{amsmath,latexsym,amsfonts,amssymb}
\usepackage{wrapfig}
\usepackage{amscd,array,epic,eepic,calc,bbm,float}
\usepackage{ifthen,psfrag,verbatim,epsfig,graphicx,enumerate,subfig}
\usepackage{makeidx}
\usepackage{amsthm}
\usepackage{mathabx,bbm,bm}
\usepackage[dvipsnames]{xcolor}
\usepackage{natbib}
\usepackage[titletoc,title]{appendix}
\usepackage{etoolbox}
\usepackage{rotating}
\usepackage{pdflscape}
\usepackage{afterpage}
\usepackage{capt-of}

\def\bql{\begin{equation}\label}
\def\eql{\end{equation}\noindent}

\def\brl{\begin{eqnarray}\label}
\def\erl{\end{eqnarray}\noindent}

\def\bro{\begin{eqnarray*}}
\def\ero{\end{eqnarray*}\noindent}

\def\brr{\begin{array}}
\def\err{\end{array}\noindent}

\def\bth{\begin{theorem}}
\def\eth{\end{theorem}}

\def\bcr{\begin{corollary}}
\def\ecr{\end{corollary}}

\def\bpr{\begin{proposition}}
\def\epr{\end{proposition}}

\def\blm{\begin{lemma}}
\def\elm{\end{lemma}}

\def\bdf{\begin{definition}}
\def\edf{\end{definition}}

\def\bas{\begin{assumptions}}
\def\eas{\end{assumptions}}

\def\bex{\begin{example}\rm}
\def\eex{\end{example}}

\def\bxx{\begin{exercise}\rm}
\def\exx{\end{exercise}}

\def\brm{\begin{remark}\rm}
\def\erm{\end{remark}}

\def\bmx{\begin{matrix}}
\def\emx{\end{matrix}}

\def\bma{\begin{pmatrix}}
\def\ema{\end{pmatrix}}

\def\bcs{\begin{cases}}
\def\ecs{\end{cases}}

\def\btb{\begin{center}\begin{tabular}}
\def\etb{\end{tabular}\end{center}}

\def\bit{\begin{itemize}}
\def\eit{\end{itemize}}

\def\bew{\par\noindent{\bf Proof }}
\def\qed{\quad\hfill$\square$}

\def\e{\epsilon}

\def\j{\psi}

\def\p{\pi}

\def\r{\rho}

\def\KC{{\cal K}}

\def\XB{{\bm X}}

\def\xb{{\bm x}}
\def\yb{{\bm y}}

\def\zerob{{\bf0}}

\def\pbb{{\mathbb P}}

\def\rbb{{\mathbb R}}

\def\nf{\infty}

\def\sm{\setminus}
\def\ss{\subset}

\def\topr{\buildrel\pbb\over\rightarrow}

\usepackage{framed}
\usepackage{comment}
\definecolor{shadecolor}{gray}{0.9}
\specialcomment{extra}{\begin{shaded}}{\end{shaded}}

\newtheorem{lemma}{Lemma}
\newtheorem{prop}{Proposition}
\newtheorem{cor}{Corollary}

\theoremstyle{definition}

\theoremstyle{definition}
\newtheorem{rmk}{Remark}
\theoremstyle{definition}
\newtheorem{example}{Example}

\newcommand{\bX}{\bm{X}}
\newcommand{\bXp}{\bm{X}_P}
\newcommand{\bXe}{\bm{X}_E}

\newcommand{\RV}{\mathrm{RV}}
\newcommand{\bbE}{\mathbb{E}}

\usepackage{authblk}

\makeindex

\newtoggle{arxiv}
\toggletrue{arxiv}

\begin{document}
\date{\today}
\title{Linking representations for multivariate extremes via a limit set}

\author[1]{Natalia Nolde}
\author[2]{Jennifer L. Wadsworth}
\affil[1]{Department of Statistics, University of British Columbia, Canada}
\affil[2]{Department of Mathematics and Statistics, Lancaster University, UK}

\renewcommand\Authands{ and }

\maketitle
\setcounter{page}{1}

\setcounter{equation}{0}

\begin{abstract}
The study of multivariate extremes is dominated by multivariate regular variation, although it is well known that this approach does not provide adequate distinction between random vectors whose components are not always simultaneously large. Various alternative dependence measures and representations have been proposed, with the most well-known being hidden regular variation and the conditional extreme value model. These varying depictions of extremal dependence arise through consideration of different parts of the multivariate domain, and particularly exploring what happens when extremes of one variable may grow at different rates to other variables. Thus far, these alternative representations have come from distinct sources and links between them are limited. In this work we elucidate many of the relevant connections through a geometrical approach. In particular, the shape of the limit set of scaled sample clouds in light-tailed margins is shown to provide a description of several different extremal dependence representations. 


\end{abstract}

\section{Introduction}

\iftoggle{arxiv}{}{}
Multivariate extreme value theory is complicated by the lack of natural ordering in $\mathbb{R}^d$, and the infinite possibilities for the underlying set of dependence structures between random variables. Some of the earliest characterizations of multivariate extremes were inspired by consideration of the vector of normalized componentwise maxima. Let $\bX =(X_1,\ldots,X_d) \in \mathbb{R}^d$ with $X_{j} \sim F_j$, and consider a sample $\bX_i =(X_{1i},\ldots,X_{di})$, $i=1,\ldots,n$, of independent copies of $\bX$. For a fixed $j$, defining $M_{j,n} = \max_{1 \leq i \leq n}(X_{ji})$, the extremal types theorem tells us that if we can find sequences such that $(M_{j,n} -b_{j,n})/a_{j,n}$ converges to a non-degenerate distribution, then this is the generalized extreme value distribution. Moreover, the sequence $b_{j,n} \sim  F_j^{-1}(1-c/n)$, $n \to \infty$, i.e., is of the same order as the $1-1/n$ quantile. A natural multivariate extension is then to examine the distribution of the vector of componentwise maxima, $(\bm{M}_n -\bm{b}_n)/\bm{a}_n$. This is intrinsically tied up with the theory of multivariate regular variation, because it leads to examination of the joint behaviour of the random vector when all components are growing at the rate determined by their $1-1/n$ quantile. If all components were marginally standardized to focus only on the dependence, then all normalizations would be the same.

In normalizing all components by the same amount, we only consider the dependence structure in a single ``direction'' in $\mathbb{R}^d$. In some cases this turns out to provide a rich description of the extremal dependence: if the limiting distribution of componentwise maxima does not have independent components, then an infinite variety of dependence structures are possible, indexed by a moment-constrained measure on a $d-1$ dimensional unit sphere. However, when the limiting dependence structure is independence, or even when some pairs are independent, this representation fails to discriminate between qualitatively different underlying dependence structures. While consideration of componentwise maxima is not necessarily a common applied methodology these days, the legacy of this approach persists: statistical methods that assume multivariate regular variation, such as multivariate generalized Pareto distributions, are still very popular in practice \citep[e.g.,][]{EngelkeHitz20}. A recent theoretical treatment of multivariate regular variation is given in \citet{KulikSoulier20}.

Various other representations for multivariate extremes have emerged that analyze the structure of the dependence when some variables are growing at different rates to others. These include the so-called conditional extreme value model \citep{HeffernanTawn04,HeffernanResnick07}, whereby the components of $\bX$ are normalized according to how they grow with a single component, $X_j$ say. Related work examines behaviour in relation to an arbitrary linear functional of $\XB$ \citep{BE07}. The conditional representation allows consideration of those regions where some or all variables grow at a lesser rate than $X_j$ if this is the region where the observations tend to lie. In other words, the limit theory is suited to giving a more detailed description of a broader range of underlying dependence structures. Another representation that explicitly considers different growth rates is that of \citet{WadsworthTawn13}. They focus particularly on characterizing joint survival probabilities under certain classes of inhomogeneous normalization; this was found to reveal additional structure that is not evident when applying a common scaling. More recently, \citet{Simpsonetal18} have examined certain types of unequal scaling with a view to classifying the strength of dependence in any sub group of variables of $\bX$.

An alternative approach to adding detail to the extremal dependence structure focuses not on different scaling orders, but rather on second order effects when applying a common scaling. This idea was introduced by \citet{LedfordTawn96}, and falls under the broader umbrella of hidden regular variation \citep{Resnick02}. Various manuscripts have focused on analogizing concepts from standard multivariate regular variation to the case of hidden regular variation \citep[e.g.,][]{RamosLedford09}, but this approach still only focuses on a restricted region of the multivariate space where all variables are large simultaneously. For this reason, although higher-dimensional analogues exist, they are often not practically useful for dimension $d>2$.

Another manner of examining the extremal behaviour of $\bm{X}$ is to consider normalizing the variables such that they converge onto a limit set \citep[e.g.,][]{Davisetal88,BalkemaNolde10}, described by a so-called gauge function. This requires light-tailed margins, which may occur naturally or through a transformation. If the margins are standardized to a common light-tailed form, then the shape of the limit set is revealing about the extremal dependence structure of the random variables, exposing in which directions we expect to see more observations.

Although various connections have been made in the literature, many of these representations remain somewhat disjointed. For example, there is no obvious connection between the conditional extremes methodology and the representation of \citet{LedfordTawn96,LedfordTawn97}, and whilst \citet{WadsworthTawn13} provided a modest connection to conditional extremes, many open questions remain. In this paper we reveal several hitherto unknown connections that can be made through the shape of the limit set and its corresponding gauge function, when it exists, and provide a step towards unifying the treatment of multivariate extremes. 

We next provide further elaboration and definition of the different representations of extremal dependence. For some definitions, it is convenient to have a standardized marginal form; we focus mainly on standard Pareto or standard exponential margins with notation $\bXp$ and $\bXe$, respectively. As mentioned above, working with common margins highlights dependence features. In Section~\ref{sec:def} we recall the formulations of various representations for multivariate extremes, and provide a thorough background to the concepts of limit sets and their gauge functions, proving a useful new result on marginalization.  
 Section~\ref{sec:connections} details connections linking conditional extremes, the representation of \citet{WadsworthTawn13}, \citet{LedfordTawn96}, and that of \citet{Simpsonetal18}. We provide \iftoggle{arxiv}{several}{}illustrative examples in Section~\ref{sec:examples} and conclude in Section~\ref{sec:Discussion}.

\section{Background and definitions}
\label{sec:def}
\subsection{Multivariate regular variation}
\label{sec:MRVintro}
A measurable function $f: \mathbb{R}_+ \to \mathbb{R}_+$ is regularly varying at infinity (respectively, zero) with index $\rho \in\mathbb{R}$ if, for any $x>0$, $f(tx)/f(t) \to x^\rho$, as $t \to \infty$  (or, respectively, $t \to 0$). We write $f \in \RV_\rho^\infty$ or $f \in \RV_\rho^0$, omitting the superscript in generic cases. If $f \in \RV_0$, then it is called slowly varying.

The random vector $\bX$ is multivariate regularly varying on the cone $\mathbb{E} = [0,\infty]^d \setminus \{\bm{0}\}$, with index $\alpha>0$, if for any relatively compact $B \subset \mathbb{E}$,
\begin{align}
 t\pbb(\bX/b(t) \in B) \to \nu(B), \qquad t \to \infty, \label{eq:mrv}
\end{align}
with $\nu(\partial B) = 0$, $b(t) \in \RV^\infty_{1/\alpha}$, and the limit measure $\nu$ homogeneous of order $-\alpha$; see, e.g., \citet{Resnick2007}, Section~6.1.4. The parts of $\mathbb{E}$ where $\nu$ places mass reveal the broad scale extremal dependence structure of $\bm{X}$. Specifically, note that we have the disjoint union $\mathbb{E} = \bigcup_{C} \mathbb{E}_C$, where
\begin{align}
 \mathbb{E}_C = \{\bm{x} \in \mathbb{E}: x_j > 0, j\in C; x_i = 0, i\not\in C\}=:(0,\infty]^C\times\{0\}^{D \setminus C}, \label{eq:EC}
\end{align}
and the union is over all possible subsets $C \subseteq D=\{1,\ldots, d\}$, excluding the empty set. If $\nu(\bbE_C)>0$ then the variables indexed by $C$ can take their most extreme values simultaneously, whilst those indexed by $D \setminus C$ are non-extreme. 

The definition of multivariate regular variation in equation~\eqref{eq:mrv} requires tail equivalence of the margins. In practice, it is rare to find variables that have regularly varying tails with common indices, and multivariate regular variation is a dependence assumption placed on standardized variables. Without loss of generality, therefore, we henceforth consider $\bX=\bXp$ with standard Pareto(1) margins in which case $\alpha=1$ and $b(t) = t$. 

Frequently, the set $B$ in \eqref{eq:mrv} is taken as $[\bm{0},\bm{x}]^c = \mathbb{E}\setminus [\bm{0},\bm{x}]$, leading to the \emph{exponent function},
\begin{align}
 V(\bm{x}) = \nu([\bm{0},\bm{x}]^c). \label{eq:V}
\end{align}
Suppose that derivatives of $V$ exist almost everywhere; this is the case for popular parametric models, such as the multivariate logistic \citep{Gumbel60}, H\"{u}sler--Reiss \citep{HueslerReiss89} and asymmetric logistic distributions \citep{Tawn90}. Let $\partial^{|C|}/\partial \bm{x}_C = \prod_{i \in C} \partial/ \partial x_i$. If the quantity $\lim_{x_j \to 0, j \not\in C}\partial^{|C|}V(\bm{x})/\partial \bm{x}_C$ is non-zero, then the group of variables indexed by $C$ places mass on $\mathbb{E}_C$ \citep{ColesTawn91}.

Multivariate regular variation is often phrased in terms of a radial-angular decomposition. If~\eqref{eq:mrv} holds, then for $r\geq1$,
\begin{align*}
\pbb(\bXp /\|\bXp\| \in A, \|\bXp\|>tr) /\pbb(\|\bXp\|>t) \to H(A)r^{-1},\qquad t\to\nf,
\end{align*}
where $A \subset \mathcal{S}=\{\bm{w}\in\mathbb{R}^d_+: \|\bm{w}\| = 1\}$ and $\|\cdot\|$ is any norm. That is, the radial variable $R= \|\bXp\|$ and the angular variable $\bm{W}=\bXp /\|\bXp\|$ are independent in the limit, with $R \sim$ Pareto(1) and $\bm{W} \in \mathcal{S}$ following distribution $H$. The support of the so-called spectral measure~$H$ can also be partitioned in a similar manner to $\mathbb{E}$. Letting
\[
  \mathbb{A}_C = \{\bm{w} \in \mathcal{S}: w_j > 0, j\in C; w_i = 0, i\not\in C\},
\]
we have $\mathcal{S} = \bigcup_C \mathbb{A}_C$. The measure $\nu$ places mass on $\mathbb{E}_C$ if and only if $H$ places mass on $\mathbb{A}_C$.

\subsection{Hidden regular variation}
Hidden regular variation arises when: (i) there is multivariate regular variation on a cone (say $\mathbb{E}$), but the mass concentrates on a subcone $\tilde{\mathbb{E}} \subset \mathbb{E}$, and (ii) there is multivariate regular variation on the subcone $\bbE' \subseteq \bbE \setminus \tilde{\mathbb{E}}$ with a scaling function of smaller order than on the full cone. Suppose that~\eqref{eq:mrv} holds, and $\nu$ concentrates on $\tilde{\mathbb{E}}$, in the sense that $\nu(\bbE \setminus \tilde{\mathbb{E}}) = 0$. For measurable $B \subset \bbE'$, we have hidden regular variation on $\bbE'$ if
\begin{align}
 t\pbb(\bXp/c(t) \in B) \to \nu'(B), \qquad t \to \infty, \qquad c(t) = o(t),\quad c(t) \in \RV_{\zeta}^\infty, \quad\zeta \in (0,1], \label{eq:hrv}
\end{align}
with $\nu'(\partial B) = 0$ and the limit measure $\nu'$ homogeneous of order $-1/\zeta$ \citep[][Section 9.4.1]{Resnick2007}.

The most common cone to consider is $\bbE' = (0,\infty]^d$. This leads to the \emph{residual tail dependence coefficient}, $\eta_D \in (0,1]$ \citep{LedfordTawn96}. That is, suppose that~\eqref{eq:hrv} holds on $(0,\infty]^d$, then the regular variation index $\zeta=\eta_D$. The residual tail dependence coefficient for the subset $C \subset D$ is found through considering cones of the form
\[
\mathbb{E}'_C = \{\bm{x} \in \mathbb{E}: x_j > 0, j\in C; x_i \in[0,\infty], i\not\in C\} =:(0,\infty]^C\times[0,\infty]^{D \setminus C},
\]
for which $\zeta= \eta_C$.

\subsection{Different scaling orders}
\label{sec:dso}

\subsubsection{Coefficients $\tau_C(\delta)$}
  \citet{Simpsonetal18} sought to examine the extremal dependence structure of a random vector through determination of the cones $\bbE_C$ for which $\nu(\bbE_C)>0$. Direct consideration of (hidden) regular variation conditions on these cones is impeded 
by the fact that $\pbb(\bXp/b(t) \in B) = 0$ for all $B \subset \bbE_C$, $C \neq D$, since no components of $\bXp/b(t)$ are equal to zero for $t< \infty$. \citet{Simpsonetal18} circumvent this issue by assuming that if $\nu(\bbE_C)>0$, then there exists $\delta<1$ such that 
\begin{align}
\lim_{t \to \infty} t\pbb(\min_{i \in C}X_{P,i}> xt, \max_{j \in D \setminus C}X_{P,j} \leq y t^\delta)=  \lim_{t \to \infty} t\pbb(\min_{i \in C}X_{P,i}/t> x, \max_{j \in D \setminus C}X_{P,j}/t \leq y t^{\delta-1})>0, \qquad x,y>0. \label{eq:SWT1}
\end{align}
As such, under normalization by $t$, components of the random vector indexed by $C$ remain positive, whereas those indexed by $D\setminus C$ concentrate at zero. Note that if assumption~\eqref{eq:SWT1} holds for some $\delta<1$ then it also holds for all $\delta' \in [\delta,1]$. \citet{Simpsonetal18} expanded assumption~\eqref{eq:SWT1} to
\begin{align}
\pbb(\min_{i \in C}X_{P,i}> xt, \max_{j \in D \setminus C}X_{P,j} \leq y t^\delta) \in \RV_{-1/\tau_C(\delta)}^{\infty}, \qquad \delta \in [0,1] \label{eq:SWT2},
\end{align}
where~\eqref{eq:SWT2} is viewed as a function of $t$, and the regular variation coefficients $\tau_C(\delta) \in (0,1]$. For a fixed $\delta$, $\tau_C(\delta)<1$ implies either that $\nu(\bbE_C)=0$, or that $\nu(\bbE_C)>0$, but that $\delta$ is too small for~\eqref{eq:SWT1} to hold; see \citet{Simpsonetal18} for further details. Considering the coefficients $\tau_C(\delta)$ over all $C$ and $\delta \in [0,1]$ provides information about the cones on which $\nu$ concentrates.

\subsubsection{Angular dependence function $\lambda(\bm{\omega})$}
\citet{WadsworthTawn13} detailed a representation for the tail of $\bXp$ where the scaling functions were of different order in each component. They focussed principally on a sequence of univariate regular variation conditions, characterizing 
\begin{align}
 \pbb(\bXp>t^{\bm{\omega}}) = \ell(t; \bm{\omega})t^{-\lambda(\bm{\omega})}, \qquad \bm{\omega} \in \mathcal{S}_{\Sigma}=\left\{\bm{\omega} \in [0,1]^d : \sum_{j=1}^d \omega_j = 1\right\}, \label{eq:WT}
\end{align}
where $\ell(t; \bm{\omega}) \in \RV_0^\infty$ for each $\bm{\omega}$ and $\lambda:\mathcal{S}_{\Sigma} \to [0,1]$. Equivalently, $\pbb(\bXe>\bm{\omega}v) = \ell(e^v; \bm{\omega})e^{-\lambda(\bm{\omega})v}$. When all components of $\bm{\omega}$ are equal to $1/d$, connection with hidden regular variation on the cone $\mathbb{E}_D$ is restored, and we have $\eta_D = d\lambda(1/d,\ldots,1/d)$. When the subcone $\bbE_D$ of $\bbE$ is charged with mass in limit~\eqref{eq:mrv}, then $\lambda(\bm{\omega}) = \max_{1\leq j \leq d}\omega_j$. One can equally focus on sub-vectors indexed by $C$ to define $\lambda_{C}(\bm{\omega})$ for $\bm{\omega}$ in a $(|C|-1)$-dimensional simplex; we continue to have $\eta_C = |C|\lambda_C(1/|C|,\ldots,1/|C|)$ and $\nu(\bbE_C)>0$ implies $\lambda_{C}(\bm{\omega}) = \max_{1 \leq j \leq |C|} \omega_j$.

\subsection{Conditional extremes}
\label{sec:CEintro}
For conditional extreme value theory \citep{HeffernanTawn04,HeffernanResnick07}, we focus on $\bXe$. Let $\bX_{E,-j}$ represent the vector $\bXe$ without the $j$th component. The basic assumption is that there exist functions $\bm{a}^j:\mathbb{R} \to \mathbb{R}^{d-1}$, $\bm{b}^j:\mathbb{R} \to \mathbb{R}_+^{d-1}$ and a non-degenerate distribution $K^j$ on $\mathbb{R}^{d-1}$ with no mass at infinity, such that
\begin{align}
 \pbb\left(\frac{\bX_{E,-j} - \bm{a}^j(X_{E,j})}{\bm{b}^j(X_{E,j})} \leq \bm{z}, X_{E,j}-t>x ~\Big|~ X_{E,j}>t\right) \to K^j(\bm{z}) e^{-x}, \qquad t \to \infty. \label{eq:ce}
\end{align}
Typically, such assumptions are made for each $j \in D$. The normalization functions satisfy some regularity conditions detailed in \citet{HeffernanResnick07}, but as \citet{HeffernanResnick07} only standardize the marginal distribution of the conditioning variable (i.e, $X_j$), allowing different margins in other variables, these conditions do not strongly characterize the functions $\bm{a}^j$ and $\bm{b}^j$ as used in~\eqref{eq:ce}. 

When joint densities exist, application of L'H\^opital's rule gives that convergence~\eqref{eq:ce} is equivalent to
\begin{align*}
\pbb\left(\frac{\bX_{E,-j} - \bm{a}^j(t)}{\bm{b}^j(t)} \leq \bm{z} ~\Big|~ X_{E,j}=t\right) \to K^j(\bm{z}), \qquad t \to \infty. 
\end{align*}
We will further assume convergence of the full joint density
\begin{align}
\frac{\partial^{d-1}}{\partial \bm{z}}\pbb\left(\frac{\bX_{E,-j} - \bm{a}^j(t)}{\bm{b}^j(t)} \leq \bm{z} ~\Big|~ X_{E,j}=t\right) \to \frac{\partial^{d-1}}{\partial \bm{z}} K^j(\bm{z}) =: k^j(\bm{z}), \qquad t \to \infty, \label{eq:cedensity}
\end{align}
which is the practical assumption needed for undertaking likelihood-based statistical inference using this model.

Connected to this approach is work in \cite{BE07}, who study asymptotic behaviour of a suitably normalized random vector~$\bX$ conditional on lying in $tH$, where $H$ is a half-space not containing the origin and $t\to\nf$. The distribution of~$\bX$ is assumed to have a light-tailed density whose level sets are homothetic, convex and have a smooth boundary. In this setting, with $H$ taken to be the vertical half-space $\{\bm{x}\in\rbb^d:\ x_d>1\}$, the limit is the so-called Gauss-exponential distribution with density $\exp\{-\bm{u}^T\bm{u}/2-v\}/(2\p)^{(d-1)/2}$, $\bm{u}\in\rbb^{d-1}$, $v>0$.

\subsection{Limit sets}
\label{sec:Gauge}
\subsubsection{Background}

\begin{figure}
	\centering
	\includegraphics[width=1\textwidth]{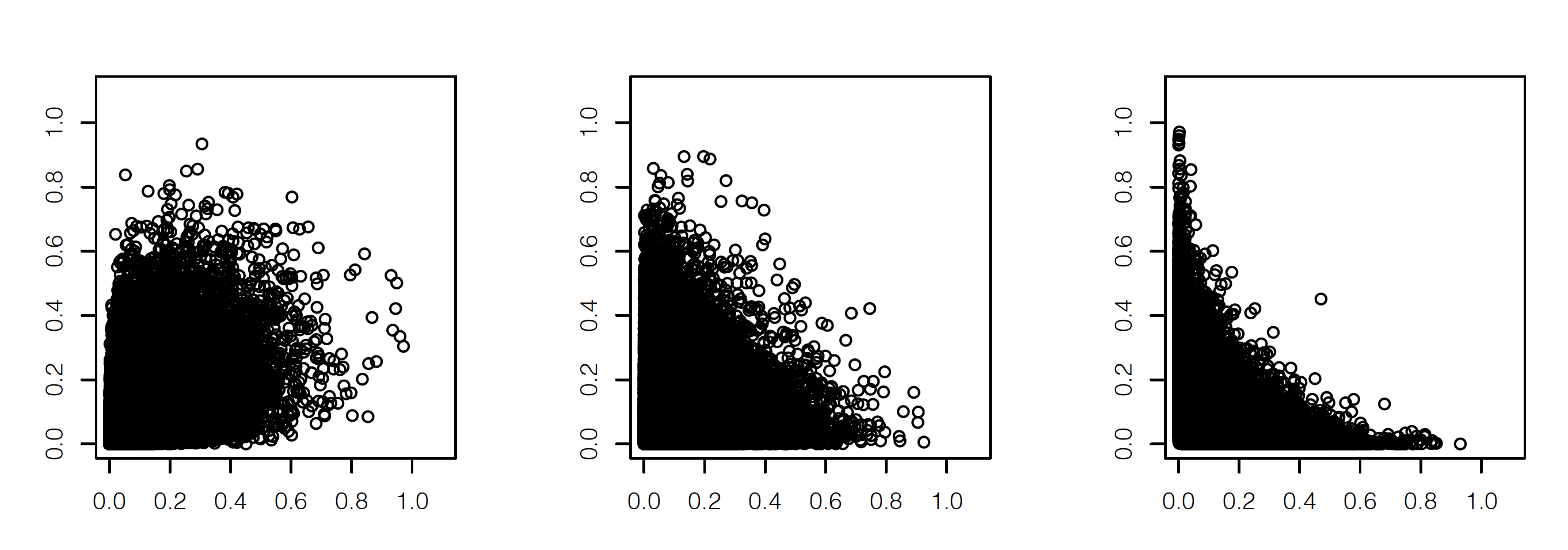}
	\caption{Sample clouds of $n=10^5$ points simulated from meta-Gaussian distributions with standard exponential margins and copula correlation parameter $\r=0.5$ (left panel), $\r=0$ (middle panel) and $\r=-0.5$ (right panel). The samples are scaled by the factor $r_n=\log n$. See Examples~\ref{exGpos} and \ref{exGneg} for details.}
	
	\label{fig:sc3}
\end{figure}

Let $\XB_1,\ldots,\XB_n$ be independent and identically distributed random vectors in $\rbb^d$. A random set $N_n=\{\XB_1/r_n,\ldots,\XB_n/r_n\}$ represents a scaled $n$-point sample cloud. We consider situations in which there exists a scaling sequence $r_n>0$, $r_n\to\nf$ such that scaled sample clouds $N_n$ converges onto a deterministic set, containing at least two points. Figure~\ref{fig:sc3} illustrates examples of sample clouds for which a limit set exists. Let $\KC_d$ denote the family of non-empty compact subsets of $\rbb^d$, and $d_H(\cdot,\cdot)$ denote the Hausdorff distance between two sets \citep{Matheron1975}. A sequence of random sets $N_n$ in $\KC_d$ converges in probability onto a \emph{limit set} $G\in\KC_d$ if $d_H(N_n,G)\topr0$ for $n\to\nf$. The following result gives convenient criteria for showing convergence in probability onto a limit set; see \citet{BEN2010}.
\begin{prop}\label{pconvLS} Random samples on $\rbb^d$ scaled by $r_n$ converge in probability onto a deterministic set $G$ in $\KC_d$ if and only if
\begin{enumerate}[(i)]
\item $n\pbb(\XB/r_n \in U^c)\to0$ for any open set $U$ containing $G$;
\item $n\pbb(\XB/r_n \in \{\bm{x}+\e B\})\to\nf$ for all $\xb\in G$ and any $\e>0$, where $B$ is the Euclidean unit ball.
\end{enumerate}
\end{prop}

Limit sets under various assumptions on the underlying distribution have been derived in \citet{Geffroy1958,Geffroy1959,Fisher1969,Davisetal88,BEN2010}. \citet{KinoshitaResnick91} give a complete characterization of the possible limit sets as well as describe the class of distribution functions for which sample clouds can be scaled to converge (almost surely) onto a limit set. Furthermore, convergence in probability onto a limit set is implied by the tail large deviation principle studied in \citet{deValk2016a,deValk2016b}.

\citet{KinoshitaResnick91} showed that if sample clouds can be scaled to converge onto a limit set almost surely, then the limit set is compact and star-shaped. A set $G$ in $\rbb^d$ is star-shaped if $\xb\in G$ implies $t\xb \in G$ for all $t \in [0,1]$. For a set $G\in\KC_d$, if the line segment $\zerob+t \xb$, $t\in [0,1)$  is contained in the interior of $G$ for every $\xb\in G$, then $G$ can be characterized by a continuous \emph{gauge function}:
$$g(\xb)=\inf\{t\ge0:\ \xb\in tG\},\qquad \xb\in\rbb^d.$$
A gauge function satisfies homogeneity: $g(t\xb)=tg(\xb)$ for all $t>0$, and the set~$G$ can be recovered from its gauge function via $G=\{\xb\in\rbb^d:\ g(\xb)\le1\}.$ Examples of a gauge function include a norm $\|\cdot\|$ on $\rbb^d$, in which case $G=\{\xb\in\rbb^d:\ \|\xb\|\le1\}$ is the unit ball in that norm.

The shape of the limit set conveys information about extremal dependence properties of the underlying distribution. In particular, \citet{BalkemaNolde10} make a connection between the shape of the limit set and asymptotic independence, whilst \citet{Nolde14} links its shape to the coefficient of residual tail dependence. We emphasize that the shape of the limit set depends on the choice of marginal distributions, as well as dependence structure. For example, if the components of $(X_1,X_2)$ are independent with common marginal distribution, then $G = \{(x,y) \in \mathbb{R}_+^2: x+y \leq 1\}$ if the margins are exponential; $G = \{(x,y) \in \mathbb{R}^2: |x|+|y| \leq 1\}$ if the margins are Laplace; and $G = \{(x,y) \in \mathbb{R}_+^2: (x^\beta +y^\beta)^{1/\beta} \leq 1\}$ if the margins are Weibull with shape $\beta>0$. In contrast, if the margins are exponential but $G$ takes the latter form, this implies some dependence between the components.

\subsubsection{Conditions for convergence onto a limit set}
Proposition~\ref{pconvLS} provides necessary and sufficient conditions for convergence onto the limit set $G$, but these conditions are not particularly helpful for determining the form of $G$ in practice.

In the following proposition, we state a criterion in terms of the joint probability density for convergence of suitably scaled random samples onto a limit set. This result is an adaptation of Proposition~3.7 in \citet{BalkemaNolde10}. The marginal tails of the underlying distribution are assumed to be asymptotically equal to a von Mises function. A function of the form $e^{-\j}$ is said to be a \emph{von Mises function} if $\j$ is a $C^2$ function with a positive derivative such that $(1/\j'(x))'\to 0$ for $x\to\nf$. This condition on the margins says that they are light-tailed and lie in the maximum domain of attraction of the Gumbel distribution, i.e., for a random sample from such a univariate distribution, coordinate-wise maxima can be normalized to converge weakly to the Gumbel distribution \citep[Proposition~1.1]{Resnick1987}.

\begin{prop}
\label{prop:logdens}
Let the random vector $\XB$ on $[0,\nf)^d$ have marginal distribution functions asymptotically equal to a von Mises function: $1-F_j(x) \sim e^{-\psi_j(x)}$ for $\psi_j(x)\sim \j(x)$, $x\to\nf$ ($j=1,\ldots,d$) and a joint probability density~$f$ satisfying: 
\bql{q1}
\dfrac{-\log f(t\xb_t)}{\j(t)}\to g_*(\xb),\qquad t\to\nf,\qquad \xb_t\to\xb,\qquad \xb\in[0,\nf)^d
\eql
for a continuous function $g_*$ on $[0,\nf)^d$, which is positive outside a bounded set. Then a sequence of scaled random samples $N_n=\{\XB_1/r_n,\ldots,\XB_n/r_n\}$ from $f$ converges in probability onto a limit set $G$ with $G=\{\xb\in[0,\nf)^d : g_*(\xb)\le 1\}$. The scaling sequence $r_n$ can be chosen as $\j(r_n)\sim \log n$. Moreover, $\max G=(1,\ldots,1)$.
\end{prop}

\bew The mean measure of $N_n$ is given by $n\pbb(\XB/r_n\in\cdot)$ with intensity $h_n(\xb) = nr_n^df(r_n\xb)$. We show the convergence of that mean measure onto $G$, implying convergence of scaled samples $N_n$; see \citet{BEN2010}, Proposition~2.3. By \eqref{q1} and the choice of~$r_n$, we have 
\begin{equation}\label{qcc}
-\log f(r_n\xb_n)/\log n\sim -\log f(r_n\xb_n)/\j(r_n)\to g_*(\xb), \quad n\to\nf,\quad \xb_n\to\xb.
\end{equation}
Continuous convergence in~\eqref{qcc} with $g_*$ continuous implies uniform convergence on compact sets. Hence, $g_*$ is bounded on compact sets. For $G=\{g_*(\xb)\le1\}$, we have $g_*(\xb)<1$ on the interior of $G$ and $g_*(\xb)>1$ on the complement of $G$.
Furthermore, applying L'H\^opital's rule and Lemma~1.2(a) in \citet{Resnick1987}, we have $$\log r_n/\j(r_n)\sim (1/\j'(r_n))/r_n\to 0,\qquad r_n\to\nf.$$
Combining these results, we see that $-\log h_n(\xb_n)\sim (g_*(\xb_n)-1)\log n$, which diverges to $-\nf$ on the interior of~$G$ and to $+\nf$ outside of~$G$. This implies that $$h_n(\xb_n)\to\bcs \nf, & \xb\in G^o,\\ 0, & \xb\in G^c,\ecs $$ 
giving convergence (in probability) of $N_n$ onto limit set~$G$.

The form of the margins $1-F_j(x) \sim e^{-\psi_j(x)}$ with $\psi_j(x) \sim \j(x)\to\nf$ gives $-\log(1-F_j(x))\sim\j(x)$; i.e., $$-\log(1-F_j(r_n))\sim \j(r_n)\sim\log n,\qquad n\to\nf.$$
This choice of $r_n$ implies that the coordinate-wise maxima scaled by $r_n$ converge in probability to 1 (\citet{G43,deHaan1970}), so that $\max G = (1,\ldots,1).$
\qed

\begin{rmk}
Condition~\eqref{q1} implies that~$-\log f$ is multivariate regularly varying on $[0,\nf)^d$. Such densities are referred to as Weibull-like. The limit function $g_*$ is homogeneous of some positive order~$k$: $g_*(t\xb)=t^k g_*(\xb)$ for all $t>0$. The gauge function $g$ of the limit set $G$ can thus be obtained from $g_*$ by setting $g(\xb)=g_*^{1/k}(\xb)$.
\end{rmk}

When the margins are standard exponential, $\j(t)=t$. Hence, for the random vector $\bXe$ with a Lebesgue density $f_E$ on $\mathbb{R}^d_+$, condition~\eqref{q1} is equivalent to 
\begin{align}\label{ass:gaugedens}
  -\log f_E(t\bm{x}_t)/t\to g_*(\bm{x}), \qquad t \to \infty,\qquad  \xb_t\to\xb,\qquad \bm{x} \in [0,\infty)^d 
 \end{align}
with the limit function $g_*$ equal to the gauge function $g$.

Whilst the assumption of a Lebesgue density might appear strict, it is a common feature in statistical practice of extreme value analysis. The assumption permits simple elucidation of the connection between different representations for multivariate extremes. Furthermore, many statistical models, including elliptical distributions and vine copulas \citep{Joe96,BedfordCooke01,BedfordCooke02}, are specified most readily in terms of their densities. 

Convergence at the density level such as in~\eqref{q1} may not always hold. The condition requires the limit function and hence the gauge function of the limit set to be continuous, excluding limit sets for which rays from the origin cross the boundary in more than one point. We provide an example of such a situation in Section~\ref{sec:examples}; see Example~\ref{exGneg}. A less restrictive set of sufficient conditions for convergence of sample clouds onto a limit set can be obtained using the survival function. The following proposition is Theorem~2.1 in~\cite{Davisetal88}, with a minor reformulation in terms of scaling.

\begin{prop}
 \label{prop:DMR}
 Suppose that the random vector $\XB$ has support on $[0,\infty)^d$, the margins are asymptotically equal to a von Mises function: $1-F_j(x) \sim e^{-\j(x)}$ for $x\to\nf$ ($j=1,\ldots,d$), and the joint survival function satisfies
 \bql{eq:DMR}
\dfrac{-\log \pbb(\XB \geq t\xb)}{\j(t)}\to g_*(\xb),\qquad t\to\nf,\qquad \xb\in [0,\infty)^d \setminus\{0\}.
\eql
Further assume that $g_*$ is strictly increasing, such that $g_*(\xb)<g_*(\yb)$ if $\xb \leq \yb$ and $\xb \neq \yb$. Then for $r_n$ satisfying $\j(r_n) \sim \log n$, the sample cloud $N_n=\{\XB_1/r_n,\ldots,\XB_n/r_n\}$ converges onto $G=\{\xb\in[0,\nf)^d : g_*(\xb)\le 1\}$.

\end{prop}

\subsubsection{Marginalization}

\label{sec:marg}

When $d>2$, a key question is the marginalization from dimension $d$ to dimension $m<d$. We prove below that, as long as the minimum over each coordinate of $g$ is well-defined, then the gauge function determining the limit set in $m$ dimensions is found through minimizing over the coordinates to be marginalized.

A continuous map $h$ from the vector space $V$ into the vector space $\widetilde V$ is \emph{positive-homogeneous} if $h(r\xb)=rh(\xb)$ for all $\xb\in V$ and all $r>0$. If $\widetilde V=\rbb^{m}$, the map $h$ is determined by the $m$ coordinate maps $h_j:V\to\rbb$, $j=1,\ldots,m$ and in this case it suffices that these maps are continuous and positive-homogeneous. 

Convergence onto a limit set is preserved under linear transformations (e.g., Lemma~4.1 in \cite{Nolde14}) and more generally under continuous positive-homogeneous maps with the same scaling sequences (Theorem~1.9 in \cite{BN20}). A consequence of the latter result, referred to as the Mapping Theorem, is that projections of sample clouds onto lower-dimensional sub-spaces also converge onto a limit set. 

\begin{prop}\label{prop:marg} Let $N_n$ be an $n$-point sample cloud from a distribution of random vector $\XB$ on $\rbb^d$. Assume $N_n$ converges in probability, as $n\to\nf$, onto a limit set $G=\{\xb\in\rbb^d : g(\xb)\le1\}$ for a gauge function $g$. Let $\widetilde\XB = (X_i)_{i\in I_m}$ denote an $m$-dimensional marginal of $\XB$, where $I_m\ss I=\{1,\ldots,d\}$ is an index set with $|I_m|=m$. Sample clouds from $\widetilde\XB$ also converge, with the same scaling, and the limit set $\tilde G = P_m(G) = \{\yb\in\rbb^m : \tilde g(\yb)\le1\}$, where $P_m$ is a projection map onto the coordinates of $\widetilde\XB$ and $$\tilde g(\yb) = \min_{\{x_i : i\in I\sm I_m\}} g(\xb),\qquad \xb=(x_1,\ldots,x_d),\qquad \yb=(x_i)_{i\in I_m}.$$
\end{prop}

\bew Consider the bivariate case first with $\widetilde\XB=X_2$. Sample clouds from $X_2$ converge onto the limit set $\tilde G\ss\rbb$, which is the projection of $G$ onto the $x_1$-coordinate axis, by the Mapping Theorem. The projection is determined by the tangent to the level curve $\{\xb\in\rbb^2 : g(\xb)=1\}$ orthogonal to the $x_1$-coordinate axis. Similarly, level curves of the gauge function $\tilde g$ of the set $\widetilde G$ are determined by tangents to the level curves $\{\xb\in\rbb^2 : g(\xb)=c\}$ for $c\in[0,1]$ orthogonal to the $x_1$-coordinate axis. These projections correspond to $x_1$ values which minimize $g(x_1,x_2)$. Sequentially minimizing over each of the coordinates to be marginalized gives the result.
\qed

An illustration of this result is given in Section~\ref{sec:egd3}.

\section{Linking representations for extremes to the limit set}
\label{sec:connections}

For simplicity of presentation, in what follows we standardize to consider exponential margins for the light-tailed case. This choice is convenient when there is positive association in the extremes, but hides structure related to negative dependence. We comment further on this case in Section~\ref{sec:Discussion}. Owing to the standardized marginals, it makes sense to refer to the limit set, rather than a limit set.

Connections between multivariate and hidden regular variation are well established, with the latter requiring the former for proper definition. Some connection between regular variation and conditional extremes was made in \citet{HeffernanResnick07} and \citet{DasResnick11}, although they did not specify to exponential-tailed margins. The shape of the limit set has been linked to the asymptotic (in)dependence structure of a random vector \citep{BalkemaNolde10,BalkemaNolde12}. Asymptotic independence is related to the position of mass from convergence~\eqref{eq:mrv} on $\bbE$, but regular variation and the existence of a limit set in suitable margins are different conditions and one need not imply the other. \citet{Nolde14} links the limit set $G$ to the coefficient of residual tail dependence, $\eta_D$.

In this section we present some new connections between the shape of the limit set, when it exists, and normalizing functions in conditional extreme value theory, the residual tail dependence coefficient, the function $\lambda(\bm{\omega})$ and the coefficients $\tau_{C}(\delta)$.

\subsection{Conditional extremes}
\label{sec:CE}

For the conditional extreme value model, the form of the normalizing functions $\bm{a}^j,\bm{b}^j$ is determined by the pairwise dependencies between $(X_{E,i},X_{E,j})$, $i \in D \setminus j$. The two-dimensional marginalization of any $d$-dimensional gauge function is given by Proposition~\ref{prop:marg}, and we simply denote this by $g$ here.

\begin{prop}
\label{prop:CEgauge}
 Suppose that for $\bXe = (X_{E,1},X_{E,2})$ convergence~\eqref{eq:cedensity}  and assumption~\eqref{ass:gaugedens} hold, where the domain of $K^j$ includes $(0,\infty)$. Define $\alpha_j = \lim_{x \to \infty} a^j(x)/x$, $j= 1,2$. Then
 \begin{enumerate}[(i)]
  \item $g(1,\alpha_1) = 1$, $g(\alpha_2, 1) = 1$.
  \item Suppose that $-\log f_E(t\bm{x}_t)/t = g(\bm{x}_t) + v(t)$, with $v(t) \in \RV_{-1}^{\infty}$ or $v(t)=o(1/t)$, and $a^j(t) = \alpha_j t + B^j(t)$ with either $B^j(t)/b^j(t) \in \RV^\infty_0$, or $B^j(t) = o(b^j(t))$. For $\beta_1,\beta_2 \leq 1$, if $g(1, \alpha_1 + \cdot) -1 \in \RV_{1/(1-\beta_1)}^0$, then  $b^1(x) \in \RV_{\beta_1}^\infty$; similarly if $g(\alpha_2 + \cdot,1) -1 \in \RV_{1/(1-\beta_2)}^0$, then $b^2(x) \in \RV_{\beta_2}^\infty$.
  \item If there are multiple values $\alpha$ satisfying $g(1,\alpha) = 1$, then $\alpha_1$ is the maximum such $\alpha$, and likewise for $\alpha_2$.
 \end{enumerate}
\end{prop}

\begin{figure}
	\centering
	\includegraphics[width=0.3\textwidth]{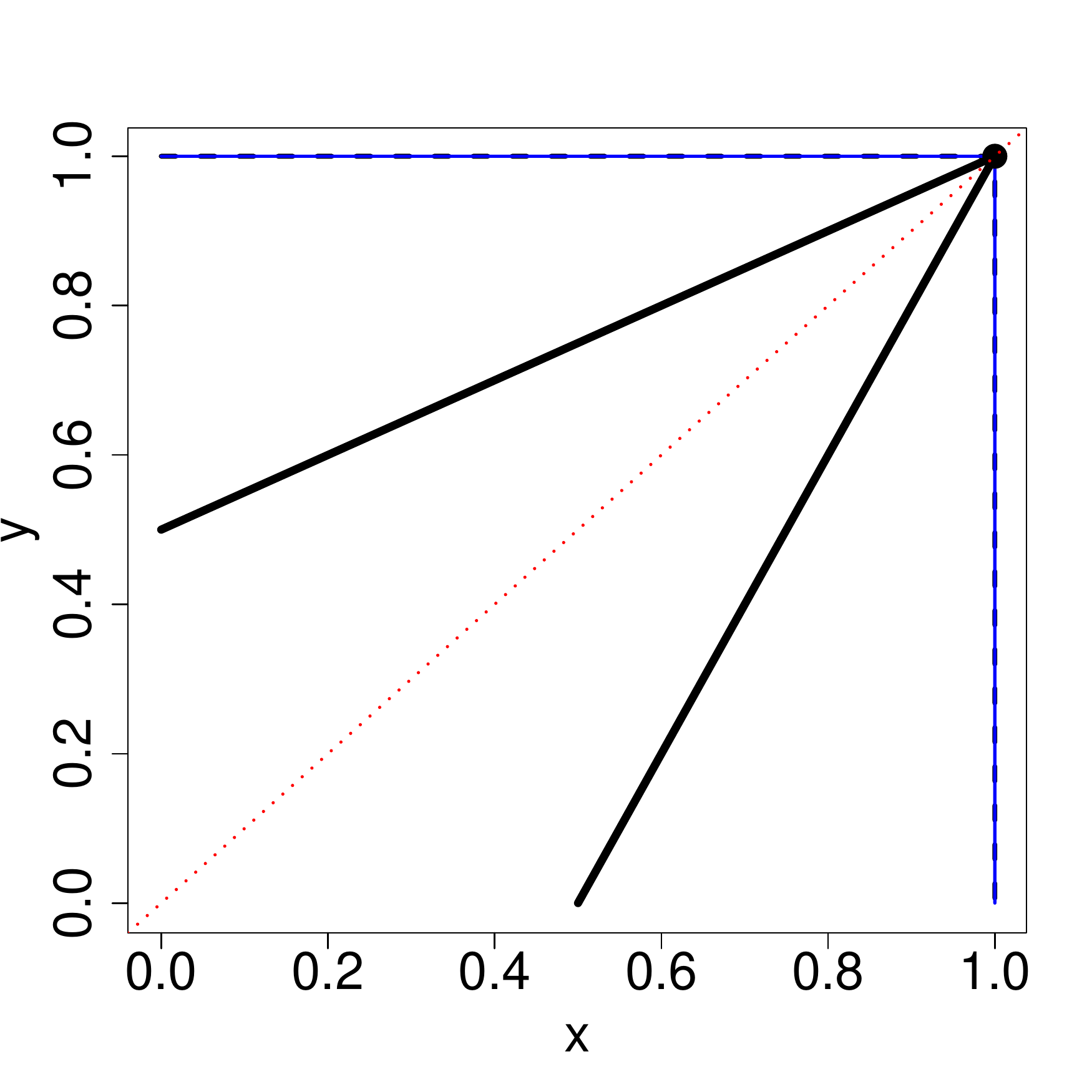}
	\includegraphics[width=0.3\textwidth]{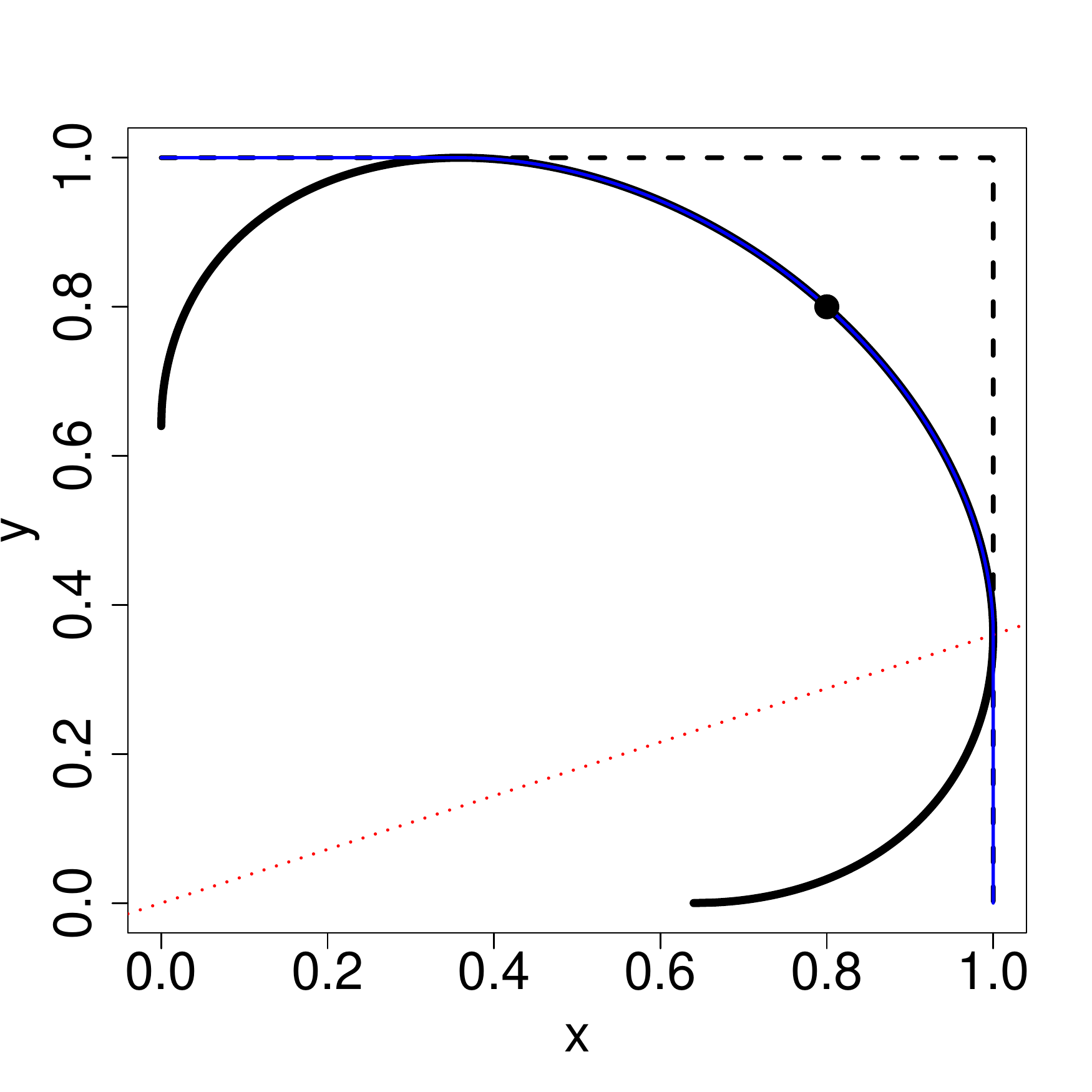}
	\includegraphics[width=0.3\textwidth]{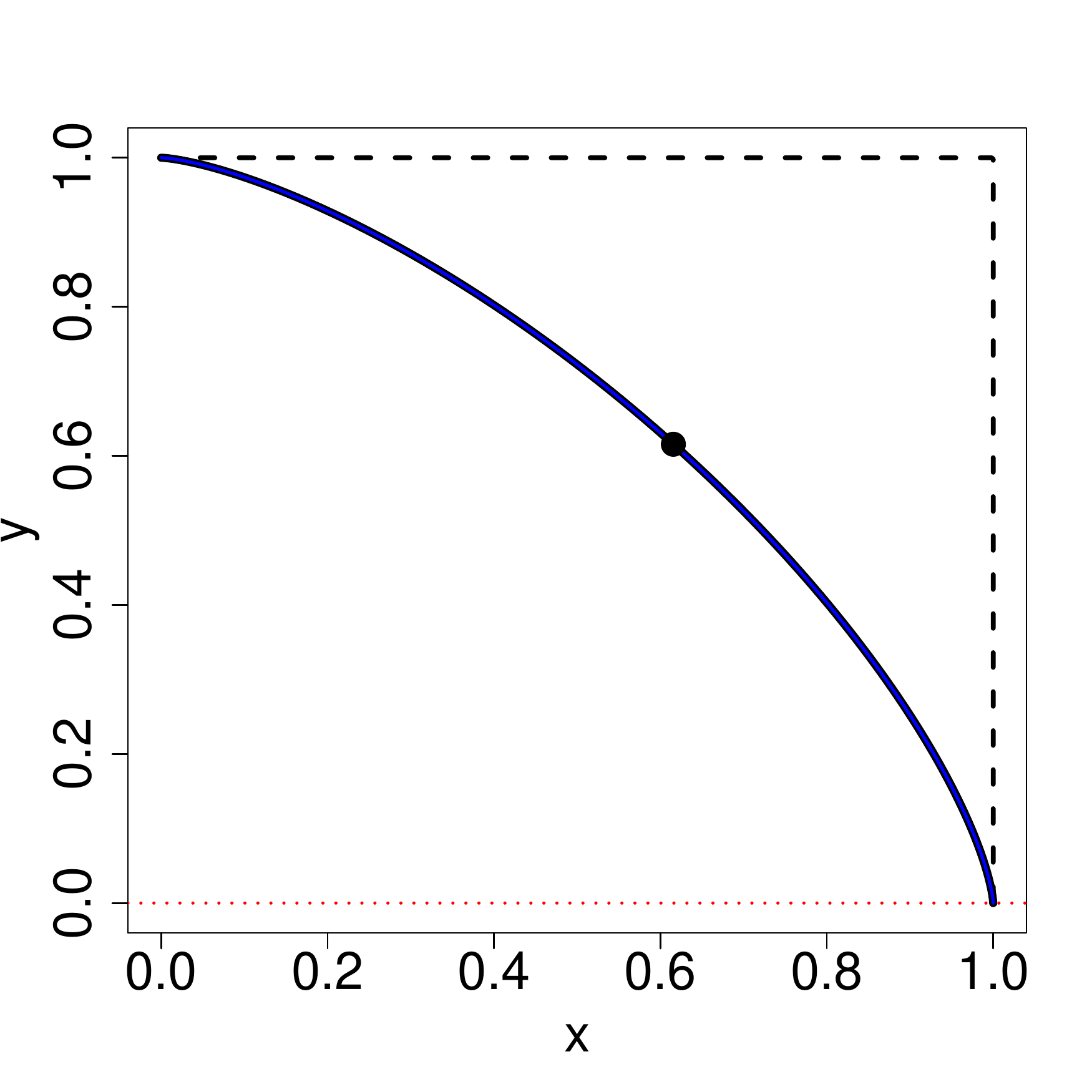}\\
	\includegraphics[width=0.3\textwidth]{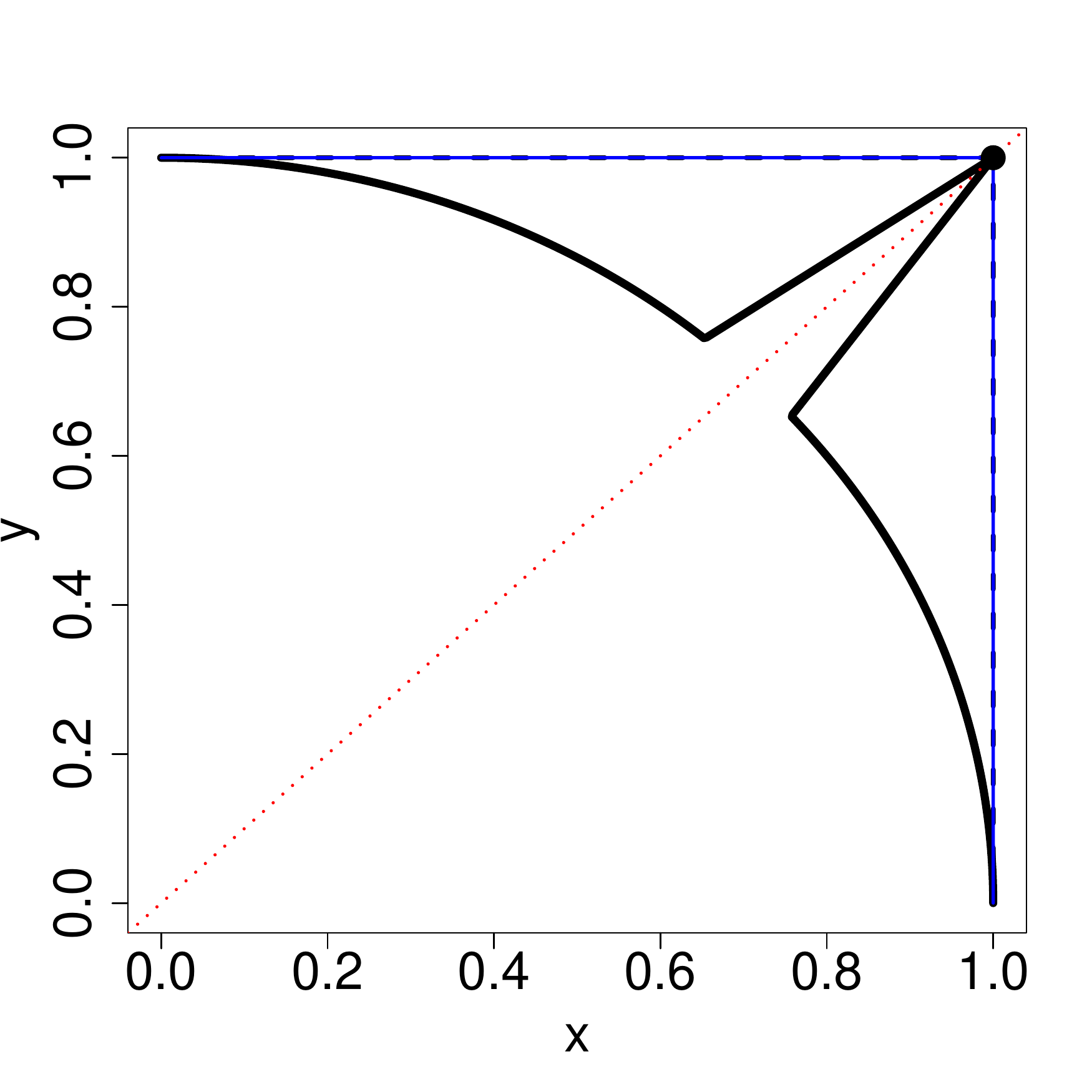}
	\includegraphics[width=0.3\textwidth]{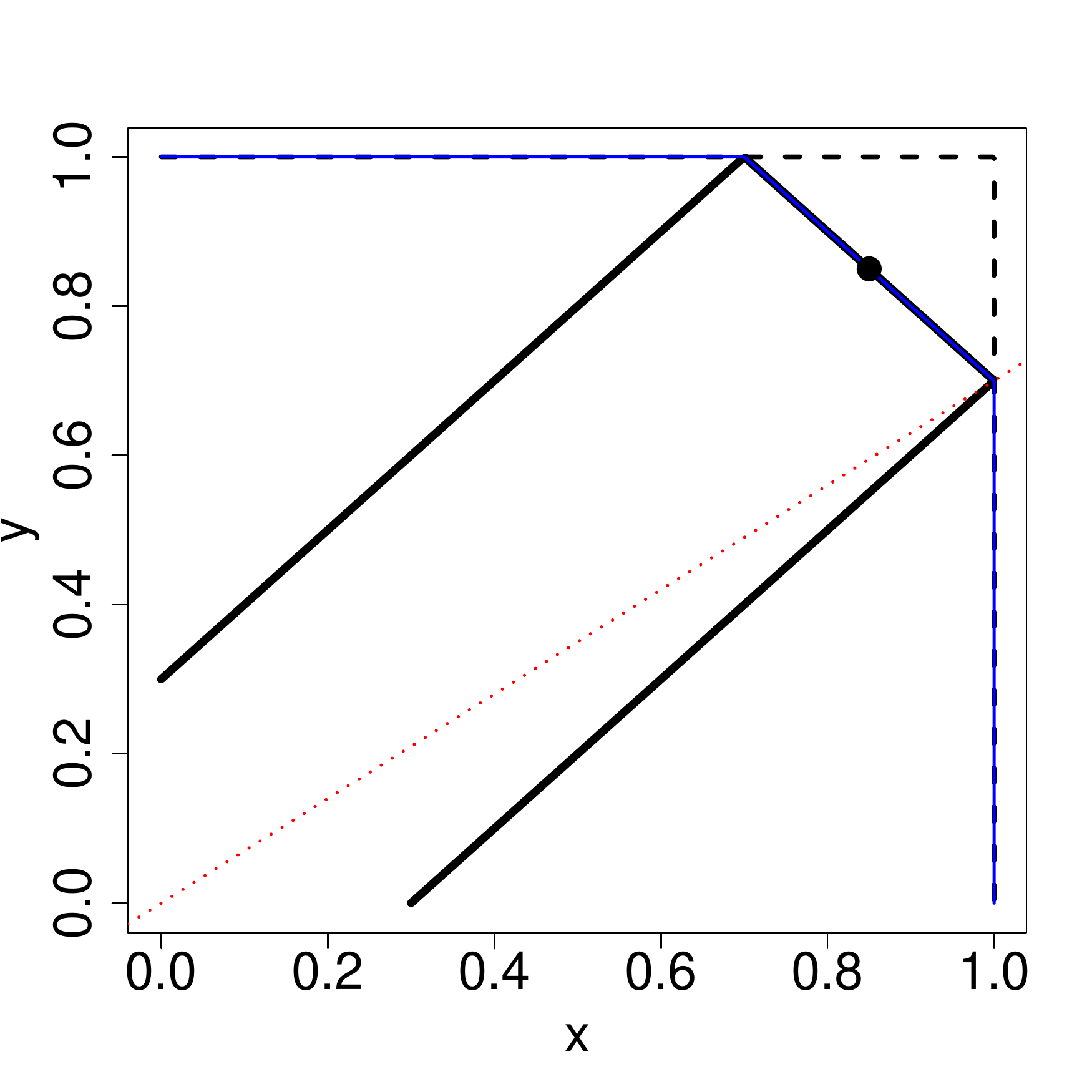}
	\includegraphics[width=0.3\textwidth]{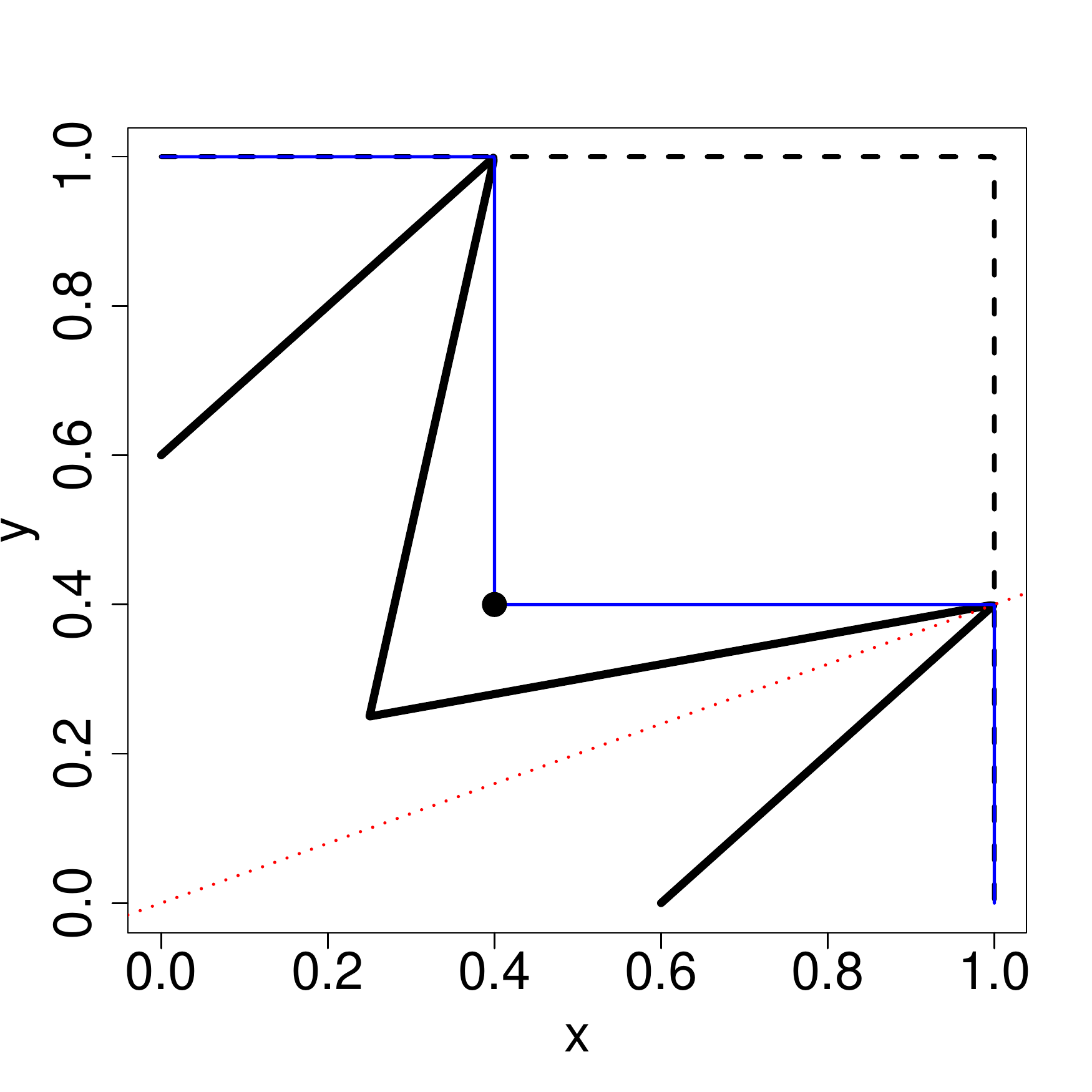}
	
	\caption{Unit level sets of six possible gauge functions for bivariate random vectors with exponential margins. The limit set $G$ is the set bounded by these level sets and the axes. In each case the red line is $y=\alpha_1 x$, the blue lines represent $\lambda(\omega,1-\omega)$ and dots $\eta_{1,2}$ (see Section~\ref{sec:lambda}). Dashed lines represent the boundary $\max(x,y)=1$. Clockwise from top left, the gauge functions represented are: (i) $\max(x,y)/\theta+(1-1/\theta)\min(x,y)$; (ii) $(x+y-2\theta\sqrt{xy})/(1-\theta^2)$; (iii) $(x^{1/\theta}+y^{1/\theta})^\theta$; (iv) $\max\{(x-y)/\theta,(y-x)/\theta,\min(x-\mu y,y-\mu x)/(1-\theta-\mu)\}$; (v) $\max((x-y)/\theta,(y-x)/\theta,(x+y)/(2-\theta))$; (vi) $\min\{\max(x,y)/\theta_1+(1-1/\theta_1)\min(x,y), (x^{1/\theta_2}+y^{1/\theta_2})^{\theta_2}\}$. In each case $\theta \in (0,1)$; in some cases the endpoints are permitted as well. For case (iv), $\theta+\mu<1$.}
	
	\label{fig:2dgauge}
\end{figure}

Before the proof of Proposition~\ref{prop:CEgauge}, we give some geometric intuition. Figure~\ref{fig:2dgauge} presents several examples of the unit level set of possible gauge functions, illustrating the shape of the limit set, for two-dimensional random vectors with exponential margins. On each figure, the slope of the red line indicates the value of $\alpha_1$; i.e., the equation of the red line is $y=\alpha_1 x$. Intuitively, conditional extreme value theory poses the question: ``given that variable $X$ is growing, how does variable $Y$ grow as a function of $X$?''. We can now see that this is neatly described by the shape of the limit set: to first order, the values of $Y$ occurring with large $X$ are determined by the direction for which $X$ is growing at its maximum rate. The necessity of a scale normalization in the conditional extreme value limit depends on the local curvature and particularly the rate at which $g(\alpha_1+u,1)$ approaches 1 as $u \to 0$. For cases (i), (iv), (v) and (vi) of Figure~\ref{fig:2dgauge}, the function approaches zero linearly in $u$: as a consequence $b^j(t) \in \RV_0^\infty$. For case (ii) the order of decay is $u^2$ and so $b^j(t) \in \RV_{1/2}^\infty$, whilst for (iii) the order is $u^{1/\theta}$ so $b^j(t) \in \RV_{1-\theta}^\infty$.

The class of distributions represented by gauge function (vi) (bottom left) can be thought of as those arising from a mixture of distributions with gauge functions (i) and (iii), up to differences in parameter values. In such an example, there are two normalizations that would lead to a non-degenerate limit in~\eqref{eq:ce}, but ruling out mass at infinity produces the unique choice $\alpha_1=\alpha_2=1$, $\beta_1=\beta_2=0$. If instead we chose to rule out mass at $-\infty$, then we would have $\alpha_1=\alpha_2 = 0$ and $\beta_1=\beta_2 = 1-\theta$.

\begin{proof}[Proof of Proposition~\ref{prop:CEgauge}] 
In all cases we just prove one statement as the other follows analogously.\\
 (i) By assumption~\eqref{ass:gaugedens}, $(X_{E,1},X_{E,2})$ have a joint density $f_E$ and so conditional extremes convergence~\eqref{eq:cedensity} can be expressed as
 \[
  b^1(t)f_E(t,b^1(t)z+a^1(t))e^t \to k^1(z) =:e^{-h^1(z)}, \qquad t \to \infty,~~~z \in [\lim_{t \to \infty}-a^1(t)/b^1(t), \infty)
 \]
with $k^1 = e^{-h^1}$ a density. Taking logs, we have 
\begin{align}
 -\log f_E(t,b^1(t)z+a^1(t)) - t -\log b^1(t) \to h^1(z),\qquad t\to\nf. \label{eq:logconv}
\end{align}
Now use assumption~\eqref{ass:gaugedens} with $\bm{x}_t =(1,x_t)= (1,a^1(t)/t + zb^1(t)/t)$. That is,
\begin{align}
 -\log f_E(t,b^1(t)z+a^1(t)) =t g(1,x)[1+o(1)] = t g(1,x_t)[1+o(1)] , \label{eq:logg}
\end{align}
with $x = \lim_{t\to\infty} a^1(t)/t + zb^1(t)/t$. As the support of $K^1$ includes $(0,\infty)$, $h^1(z)<\infty$ for all $z \in (0,\infty)$, and combining~\eqref{eq:logconv} and~\eqref{eq:logg} we have
\begin{align}
g(1,x_t)[1+o(1)] = 1 + h^1(z)/t +\log b^1(t)/t + o(1/t).\label{eq:prelim}
\end{align}
Suppose that $b^1(t)/t \to \gamma>0$. Then $x_t \to \alpha_1 + \gamma z$, and taking $t\to \infty$ in~\eqref{eq:prelim} leads to $g(1,\alpha_1 + \gamma z) =1$ for any $z$. But since the coordinatewise supremum of $G$ is $(1,1)$, $g(x,y) \geq \max(x,y)$ which would entail $z \leq (1-\alpha_1)/\gamma$. No such upper bound applies, so we conclude $\gamma=0$, i.e., $b^1(t) = o(t)$. Now taking limits in~\eqref{eq:prelim} leads to $g(1,\alpha_1)=1.$

(ii) Let $g(1,\alpha_1+u)-1 =:r(u) \in \RV_\rho^0$, $\rho>0$. We also have from~\eqref{eq:prelim} $g(1,\alpha_1 + b^1(t)/t+B^1(t)/t)-1 =  h^1(1)/t + \log b^1(t)/t - v(t) + o(1/t)$, so that the function $b^1(t)$ is a solution to the equation
\begin{align}
 r(b^1(t)/t+B^1(t)/t) = h^1(1)/t + \log b^1(t)/t - v(t)+ o(1/t). \label{eq:exprv}
\end{align}
Equation~\eqref{eq:exprv} admits a solution if $b^1$ is regularly varying at infinity. A rearrangement provides that
\begin{align*}
 b^1(t) = t r^{-1}(h^1(1)/t + \log b^1(t)/t - v(t)+ o(1/t))[1+ B^1(t)/b^1(t)]^{-1};
\end{align*}
if $b^1$ is regularly varying then $\log b^1(t)/t \in \RV_{-1}^{\infty}$, so that using the fact that $v(t) \in \RV_{-1}^\infty$ or $v(t)=o(1/t)$, combined with $r^{-1} \in \RV_{1/\rho}^0$, yields $b^1(t) \in \RV_{1-1/\rho}^{\infty}$. We now argue that such a solution is unique in this context. We know that the normalization functions $a^1,b^1$ lead to a non-degenerate distribution $K^1$ that places no mass at infinity. By the convergence to types theorem (\citet[][p.7]{Leadbetter83}, see also part (iii) of this proof), any other function $\tilde{b}^1$ leading to a non-degenerate limit with no mass at infinity must satisfy $\tilde{b}^1(t) \sim d b^1(t)$, $t \to \infty$, for some $d>0$, so that $\tilde{b}^1 \in \RV_{1-1/\rho}^\infty$ also. Finally, setting $\beta_1=1-1/\rho$ gives $b^1 \in \RV_{\beta_1}^{\infty}$.

(iii) Suppose that
 \[
\pbb\left(\frac{X_{E,2}-a^1(t)}{b^1(t)} \leq z \big| X_{E,1}>t\right) \to K^1(z), ~~~~ \pbb\left(\frac{X_{E,2}-\tilde{a}^1(t)}{\tilde{b}^1(t)} \leq z \big| X_{E,1}>t\right) \to \tilde{K}^1(z),
 \]
where neither $K^1$ nor $\tilde{K}^1$ has mass at $+ \infty$. Then by the convergence to types theorem, $\tilde{a}^1(t) = a^1(t) + cb^1(t) +o(b^1(t))$ and $\tilde{b}^1(t) = db^1(t) +o(b^1(t))$, for some $d>0$, and $\tilde{K}^1(z) = K^1(z/d+c)$. As such, $\tilde{a}^1(t)/t \sim a^1(t)/t \sim \alpha_1$. We conclude that if there was a non-degenerate $\tilde{K}^1$ limit for which $\tilde{a}^1(t)/t \sim \tilde{\alpha}_1>\alpha_1$ then $K^1$ must place mass at $+\infty$; since by assumption it does not, then $\alpha_1$ is the maximum value satisfying $g(1,\alpha_1) = 1$.
\end{proof}

For distributions whose sample clouds converge onto a limit set described by a gauge function with piecewise-continuous partial derivatives possessing finite left and right limits, further detail can be given about $\beta_j$.

\begin{prop}
	Let $G$ be a limit set whose gauge function $g$ has piecewise-continuous partial derivatives $g_1(x,y)=\partial g(x,y)/\partial x$, $g_2(x,y)=\partial g(x,y)/\partial y$ possessing finite left and right limits, and for which the conditions of Proposition~\ref{prop:CEgauge} hold. Then (i) $\beta_1\geq0$ if $g_2(1,(\alpha_1)_+)=0$; (ii) $\beta_1=0$ if $0<g_2(1,(\alpha_1)_+)<\infty$. Further, if $\alpha_1>0$ then  $0 \leq g_2(1,(\alpha_1)_+)< \infty$, so that $\beta_1 \geq 0$. Analogous statements hold for $\alpha_2,\beta_2$.
	\label{prop:betapos}
\end{prop}

\begin{proof}
	Consider the partial derivative
	\[
	g_2(1,(\alpha_1)_+) = \lim_{u\to 0^+} \frac{g(1,\alpha_1+u)-g(1,\alpha_1)}{u} \geq 0,
	\]
	as $g(1,\alpha_1+u) \geq g(1,\alpha_1)$. We note $g(1,\alpha_1)=1$, such that $g(1,\alpha_1+u)-1 \sim u g_2(1,(\alpha_1)_+)$, $u \to 0^+$. Since this is regularly varying with index $1/(1-\beta_1)$ by assumption, $g_2(1,(\alpha_1)_+)=0$ implies $g(1,\alpha_1+u)-1=o(u)$, hence $1/(1-\beta_1)\geq1$, and $0<g_2(1,(\alpha_1)_+)<\infty$ implies $1/(1-\beta_1)=1$, so (i) and (ii) follow. If $g$ is differentiable at the point $(1,\alpha_1)$, then since $g(1,y)\geq 1$, $g_2(1,(\alpha_1)_+) = g_2(1,(\alpha_1)_-)=0$ and (i) holds. Otherwise, in a neighbourhood of $(1,\alpha_1)$, we can express
	\begin{align*}
	g(x,y) = \begin{cases}
	\hat{g}(x,y),& y \leq \alpha_1 x\\
	\tilde{g}(x,y),& y \geq \alpha_1 x,
	\end{cases}
	\end{align*}
	where the homogeneous functions $\tilde{g}$ and $\hat{g}$ have continuous partial derivatives at $(1,\alpha_1)$.
	 Euler's homogeneous function theorem gives $1 = \tilde{g}_1(1,\alpha_1)+\tilde{g}_2(1,\alpha_1)\alpha_1 = g_1(1_-,\alpha_1)+\alpha_1 g_2(1,(\alpha_1)_+)$ so that for $\alpha_1>0$, $g_2(1,(\alpha_1)_+)<\infty$, and hence (i) or (ii) hold.
	\end{proof}

We remark on links with existing work on conditional extreme value limits for variables with a polar-type representation, whereby $(X_1,X_2)=R(W_1,W_2)$ for $R>0$ and $(W_1,W_2)$ constrained by some functional dependence. \citet{AbdousEtAl05}, \citet{FougeresSoulier10} and \citet{Seifert14} consider a type of conditional extremes limit for certain such polar constructions, where in the light-tailed case, the shape of the constraint on $(W_1,W_2)$ feeds into the normalization and limit distribution. However, limit sets are sensitive to marginal choice, and because the above papers do not consider conditional extreme value limits in standardized exponential-tailed margins, further connections are limited.

\subsection{Different scaling orders: $\lambda(\bm{\omega})$}
\label{sec:lambda}
We now focus on the connection with $\lambda(\bm{\omega})$, as defined in Section~\ref{sec:dso}. When $\bm{\omega} = (1/d,\ldots,1/d)$, this yields the link with the residual tail dependence coefficient $\eta_{D}$, which has already been considered in \citet{Nolde14}. Define the region 
\[
  R_{\bm{\omega}}=\left(\frac{\omega_1}{\max(\omega_1,\ldots,\omega_d)} , \infty\right] \times \cdots \times \left(\frac{\omega_d}{\max(\omega_1,\ldots,\omega_d)}, \infty\right].
\]

\begin{prop}
	Suppose that the sample cloud $N_n = \{\bXe^1/\log n,\ldots,\bXe^n/\log n\}$ converges onto a limit set $G$, and that for each $\bm{\omega}\in\mathcal{S}_{\Sigma}$, equation~\eqref{eq:WT} holds. Then 
\label{prop:lambda}
 \[
 \lambda(\bm{\omega}) = \max(\bm{\omega}) \times r_{\bm{\omega}}^{-1} ,
 \]
 where
 \[
 r_{\bm{\omega}}=\min\left\{r \in [0,1] : r R_{\bm{\omega}} \cap G = \emptyset \right\}.
\]
\end{prop}
\begin{cor}[\citet{Nolde14}]
 \[
1/\eta_{D} = d \lambda(1/d,\ldots,1/d) =r_{(1/d,\ldots,1/d)}^{-1} = \left[\min\left\{r \in [0,1] :  \left(r , \infty\right]^d \cap G = \emptyset \right\}\right]^{-1}.
 \]
\end{cor}

\begin{proof}[Proof of Proposition~\ref{prop:lambda}]
 The proof follows very similar lines to Proposition~2.1 of \citet{Nolde14}. Firstly note that $\lambda(\bm{\omega}) = \kappa(\bm{\omega})$, where $\kappa: [0,\infty)^d \setminus \{\bm{0}\} \to (0,\infty)$ is a 1-homogeneous function defined by
 \[
   \pbb(\bX_{E}> \bm{\beta} t) = \ell(e^{t}; \bm{\beta})e^{-t \kappa(\bm{\beta})},\qquad \ell(\cdot;\bm{\beta}) \in \RV_{0}^\infty~\text{for every}~ \bm{\beta} \in [0,\infty)^d \setminus \{\bm{0}\} .
 \]
As a consequence,
 \begin{align}
  \frac{\lambda(\bm{\omega})}{\max(\bm{\omega})} = \lim_{t\to\infty}-\log\pbb\{\bX_{E} > t\bm{\omega}/\max(\bm{\omega}) \}/t. \label{eq:lambdalim}
 \end{align}
Without loss of generality, suppose that $\max(\bm{\omega}) = \omega_d$, so that $R_{\bm{\omega}}=(\omega_1/\omega_d,\infty] \times \cdots \times (\omega_{d-1}/\omega_d,\infty] \times (1,\infty]$. Because of the convergence of the sample cloud onto $G$, we have by Proposition~\ref{pconvLS} that for any $\epsilon >0$ and large enough $t$
\[
 \pbb(\bXe \in t e^\epsilon r_{\bm{\omega}} R_{\bm{\omega}}) \leq \pbb(\bXe \in t R_{(0,\ldots,0,1)}) = e^{-t} \leq \pbb(\bXe \in t e^{-\epsilon} r_{\bm{\omega}} R_{\bm{\omega}}),
\]
implying $-\log\pbb(\bXe \in t r_{\bm{\omega}} R_{\bm{\omega}}) \sim t$. Therefore $-\log\pbb(\bXe \in t R_{\bm{\omega}}) \sim t r_{\bm{\omega}}^{-1}$, and combining with equation~\eqref{eq:lambdalim} gives the result.
 
\end{proof}

\begin{figure}[h]
\centering
\includegraphics[width=0.4\textwidth]{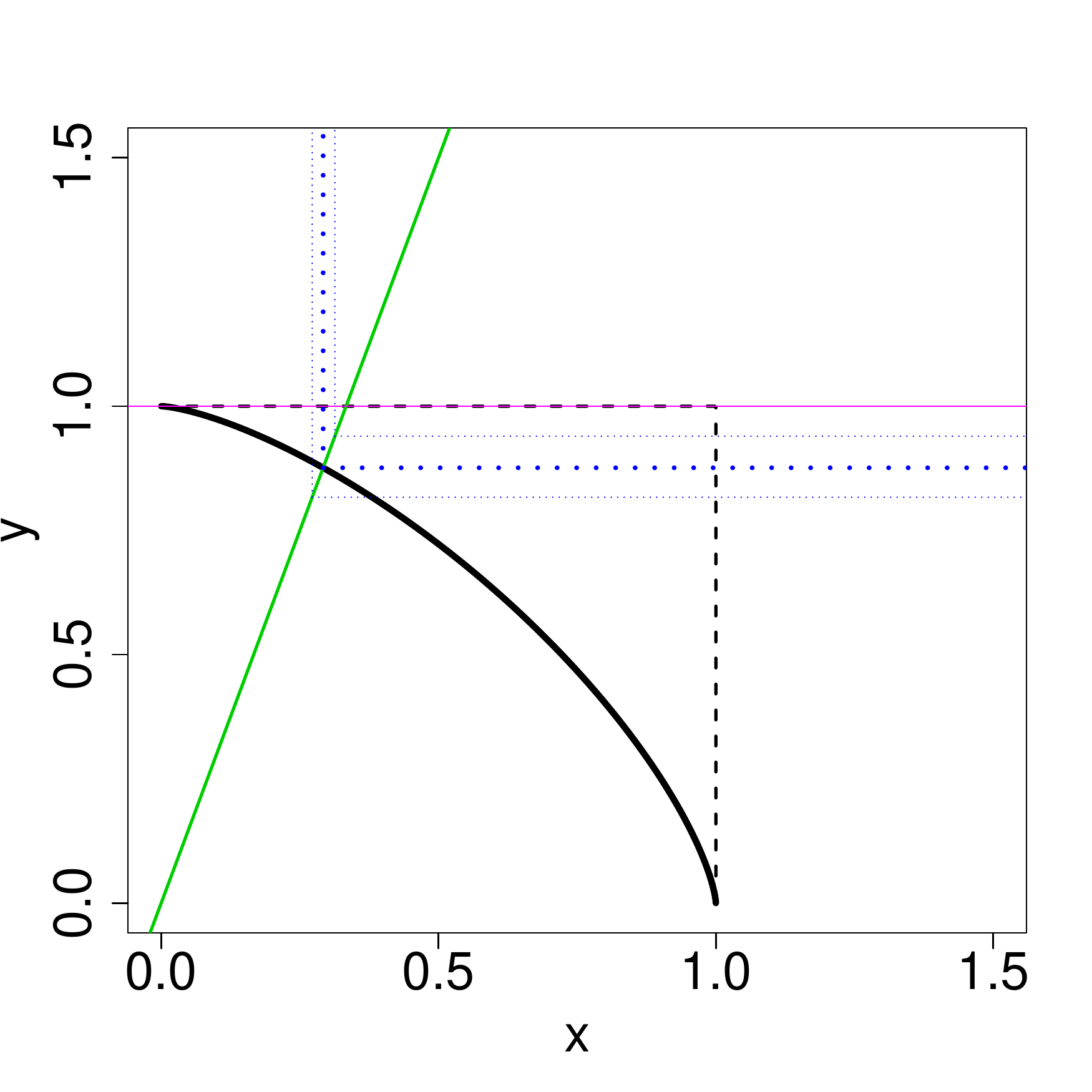}
\caption{Illustration of the concepts used in the proof of Proposition~\ref{prop:lambda}:  the green line represents the ray $x=\{\omega/(1-\omega)\}y$, $\omega<1/2$. The region above the purple line represents $R_{(0,1)}$, the region to the north-east of the thick dotted blue lines represents $r_{\omega}R_{\omega}$, whilst the two sets of thin dotted blue lines illustrate the regions $e^{-\epsilon }r_{\omega}R_{\omega}$ and $e^{\epsilon}r_{\omega}R_{\omega}$. The ratio of the distance from the origin to where the green line intersects the boundary of the limit set $G$, and where the green line intersects the boundary $\max(x,y)=1$, is equal to $r_{\omega}$.}
\label{fig:lambda}
\end{figure}

Figure~\ref{fig:lambda} illustrates some of the concepts used in the proof of Proposition~\ref{prop:lambda} when $d=2$ and $\bm{\omega} = (\omega,1-\omega)$.

The blue lines in Figure~\ref{fig:2dgauge} represent $\lambda(\bm{\omega})$, depicting the unit level set of $\lambda(\omega,1-\omega)/ \max(\omega,1-\omega)$, and the dots illustrate the value of $r_{1/2} = \eta_{1,2}$. We can now see clearly how, in two dimensions, different dependence features are picked out by the conditional extremes representation and hidden regular variation based on $\eta_{1,2}$. Often, values of $\eta_{1,2}>1/2$ or $\alpha>0$ are associated with positive extremal dependence. From example~(iv) of Figure~\ref{fig:2dgauge} (bottom right), we observe $\eta_{1,2}<1/2$ but $\alpha>0$. We have that $Y$ does grow with $X$ (and vice versa) but only at a specific rate. On the other hand, joint extremes, where $(X,Y)$ take similar values, are rare, occurring less frequently than under independence.

From example~(iv) we can also see that one of the conclusions following Proposition~2.1 in \citet{Nolde14} is not true: the point $(r_{1/2},r_{1/2})$ need not lie on the boundary of $G$, meaning that we do not necessarily have $\eta_{D} = 1/g(\bm{1})$, although we can deduce the bound $\eta_{D} \geq 1/g(\bm{1})$. Similarly, there are occasions when $g(r_{\bm{\omega}}\bm{\omega}/\max(\bm{\omega})) = 1$, implying $\lambda(\bm{\omega})=g(\bm{\omega})$, but clearly this is not always true. In Proposition~\ref{prop:lambdag}, we resolve when this is the case by representing $r_{\bm{\omega}}$ in terms of $g$. 

Define $B_{\bm{\omega}}$ to be the boundary of the region $R_{\bm{\omega}}$, i.e.,
\[
 B_{\bm{\omega}} = \bigcup_{i=1}^d \{\bm{x} \in \mathbb{R}^d_+: x_i = \omega_i/\max(\bm{\omega}), x_j \geq \omega_j/\max(\bm{\omega}), j \neq i\}.
\]

\begin{prop}
	Assume the conditions of Proposition~\ref{prop:lambda}. Then
\label{prop:lambdag}
 \[r_{\bm{\omega}} = \left[\min_{\bm{y} \in B_{\bm{\omega}}} g(\bm{y})\right]^{-1},~~~\text{and hence}~~~ \lambda(\bm{\omega}) = \max(\bm{\omega}) \times \min_{\bm{y} \in B_{\bm{\omega}}} g(\bm{y}).
 \]
\end{prop}

From Proposition~\ref{prop:lambdag}, we observe that $\lambda(\bm{\omega}) = g(\bm{\omega})$ if $\arg\min_{y \in B_{\bm{\omega}}} g(\bm{y}) = \bm{\omega}/\max(\bm{\omega})$, i.e., the vertex of the set $B_{\bm{\omega}}$. The proof of Proposition~\ref{prop:lambdag} is deferred until after Proposition~\ref{prop:taug}, for which the proof is very similar.

\begin{rmk}
	We note that $\min_{\bm{y} \in B_{\bm{\omega}}} g(\bm{y}) = \min_{\bm{y} \in R_{\bm{\omega}}} g(\bm{y})$.
\end{rmk}

\subsection{Coefficients $\tau_C(\delta)$}
\label{sec:tau}
\subsubsection{Connections to limit set $G$}
\label{sec:tau2}

\begin{figure}
\centering
 \includegraphics[width=0.3\textwidth]{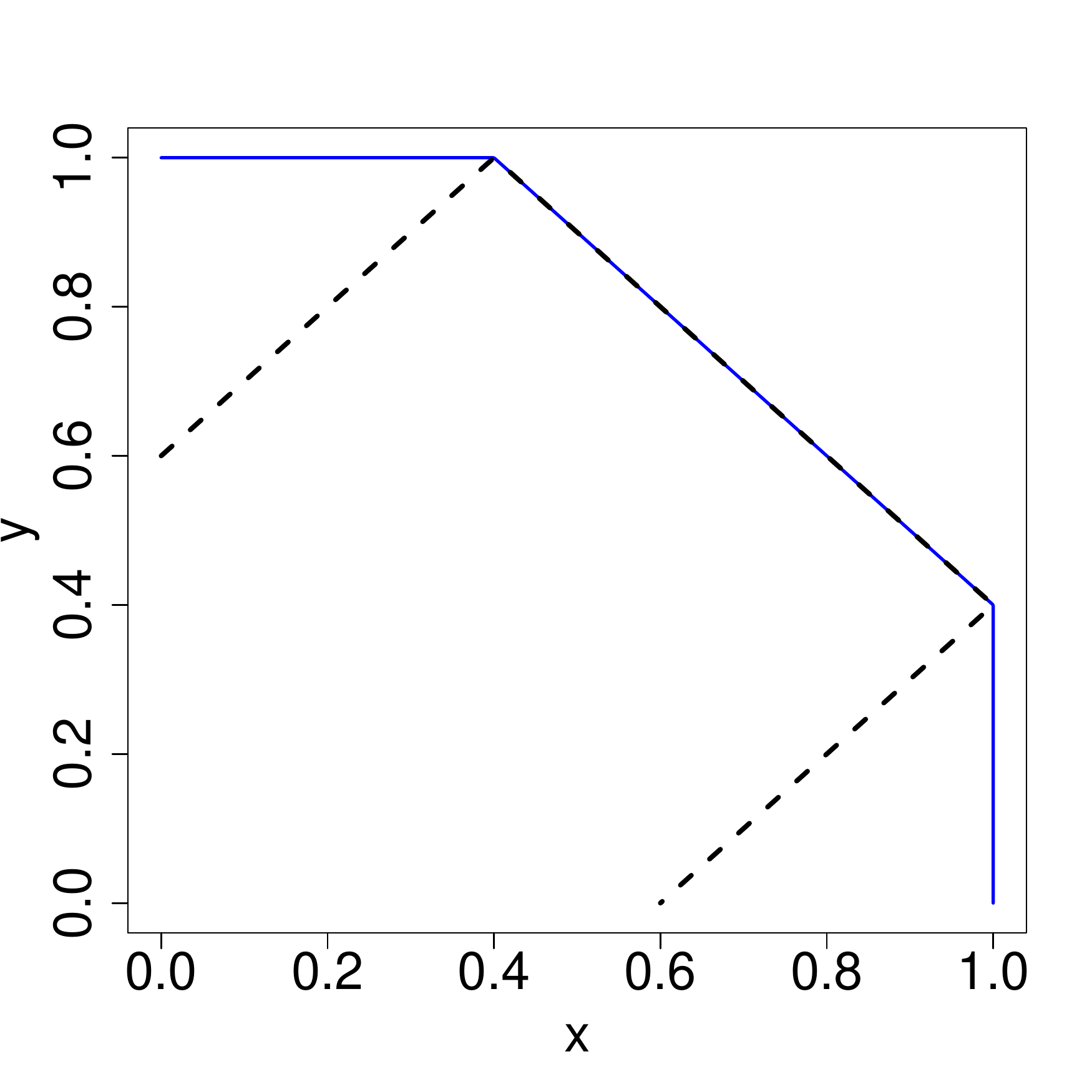}
 \includegraphics[width=0.3\textwidth]{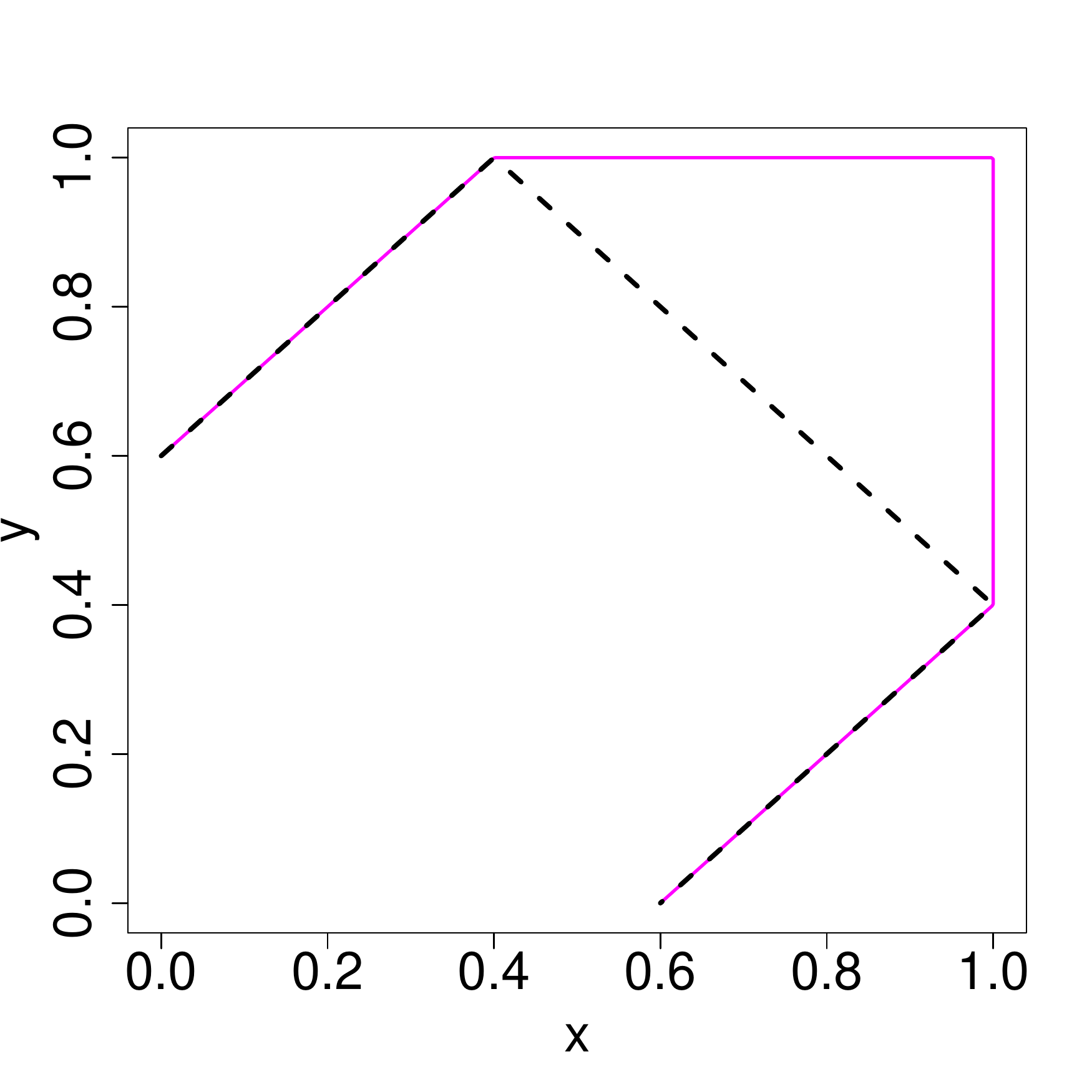}
 \includegraphics[width=0.3\textwidth]{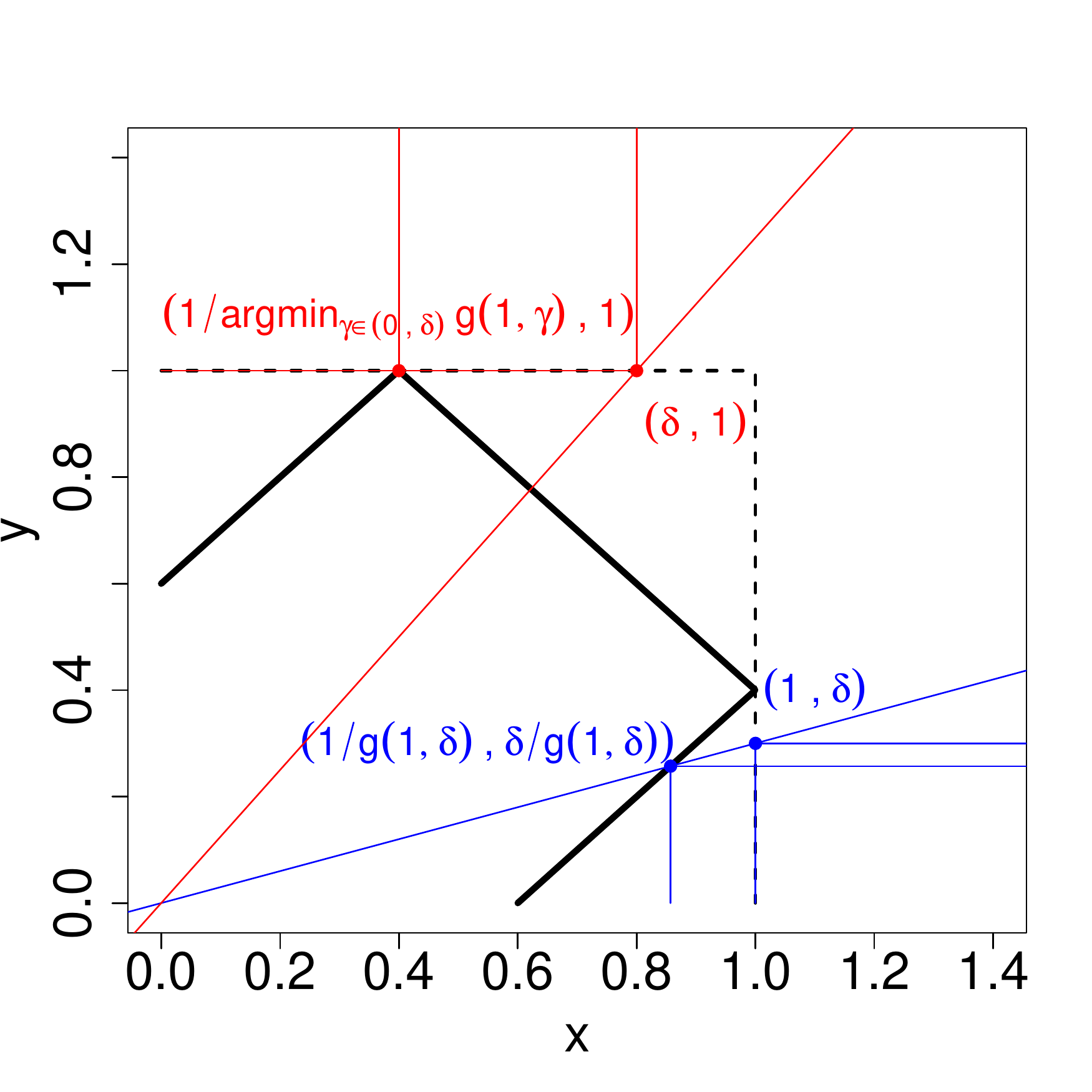}
 \caption{Illustration of $\lambda(\bm{\omega})$ and $\tau_{j}(\delta)$ for the gauge function $g(x,y) = \max\{(x-y)/\theta,(y-x)/\theta,(x+y)/(2-\theta)\}$. The left panel illustrates $\lambda(\omega,1-\omega)$ in blue. The centre panel illustrates $\tau_j(\delta)$ in purple: $\tau_1(\delta), \delta \in [0,1]$ is represented by the values below the main diagonal, whilst $\tau_2(\delta)$ is represented by those values above. The set $G$ is added on both panels with dashed lines. The right panel illustrates $\tau_1(\delta)$ and $\tau_2(\delta)$ in terms of the gauge function.}
 \label{fig:lambdatau}
\end{figure}

In two dimensions, the coefficients $\tau_{1}(\delta)$ and $\tau_{2}(\delta)$ provide a somewhat complementary concept to the function $\lambda(\bm{\omega})$. Rather than considering the impact of the limit set $G$ on the shape of the function defined by both variables exceeding thresholds growing at different rates, we are considering what is occurring when one variable exceeds a growing threshold and the other is upper bounded by a certain lesser growth rate. The left and centre panels in Figure~\ref{fig:lambdatau} provide an illustration of $\lambda(\bm{\omega})$ and $\tau_j(\delta)$ in two dimensions.

Define the region $R_{C,\delta} = (1,\infty]^C \times [0,\delta]^{D \setminus C} = \{\bm{x}: x_i \in (1,\infty], i \in C, x_j \in [0,\delta], j \in D\setminus C\}$, so that, for example, when $d=3$, $R_{\{1,3\},\delta} = (1,\infty] \times [0,\delta] \times (1,\infty]$.
 
\begin{prop}
\label{prop:tau}
Suppose that the sample cloud $N_n = \{\bXe^1/\log n,\ldots,\bXe^n/\log n\}$ converges onto a limit set $G$, and that the assumption in equation~\eqref{eq:SWT2} holds. For $\delta \in [0,1]$, and $C \subset D$,
 \begin{align*}
\tau_C(\delta) &= r_{C,\delta} = \min\left\{r \in [0,1] : r R_{C,\delta} \cap G =\emptyset \right\}.
 \end{align*}
 The coefficient $\tau_{D} = \eta_{D}$, and does not depend on $\delta$.
 \end{prop}
 
\begin{proof}
The coefficient $\tau_{D}$ describes the order of hidden regular variation on the cone $(0,\infty]^d$, which is precisely the same as $\eta_{D}$. For $\tau_C(\delta)$, $|C|<d$, we consider the function of $t$
\begin{align*}
\pbb(\min_{i \in C} X_{P,i}>tx, \max_{j \in D \setminus C} X_{P,j} \leq y t^\delta) \in \RV_{-1/\tau_C(\delta)}^\infty, \qquad 0<x,y<\infty.
\end{align*}
Take $x=y=1$. Then
\begin{align}
 \tau_C(\delta) &= \lim_{t \to \infty} \frac{-\log\pbb(X_{P,1}>t)}{-\log \pbb(\min_{i \in C} X_{P,i}>t, \max_{j \in D \setminus C} X_{P,j} \leq  t^\delta)} \notag \\&= \lim_{t \to \infty} \frac{-\log\pbb(X_{E,1}>t)}{-\log\pbb(\min_{i \in C} X_{E,i}>t, \max_{j \in D \setminus C} X_{E,j} \leq \delta t)} ,
\label{eq:tau}
 \end{align}
where the denominator in~\eqref{eq:tau} can be expressed $-\log\pbb(\bm{X}_E \in t R_{\delta})$. As in the proof of Proposition~\ref{prop:lambda}, the convergence onto the limit set and exponential margins enables us to conclude that $-\log \pbb(\bm{X}_E \in t r_{C,\delta} R_{\delta}) \sim t$, and hence  $-\log \pbb(\bm{X}_E \in t R_{\delta}) \sim tr_{C,\delta}^{-1}$. Combining with~\eqref{eq:tau} gives $\tau_C(\delta) = r_{C,\delta}$.
\end{proof}

In the two-dimensional case, it is possible to express $\tau_j(\delta)$ simply in terms of the gauge function. For higher dimensions, we refer to Proposition~\ref{prop:taug}.
\begin{prop}
	Assume the conditions of Proposition~\ref{prop:tau}. When $d=2$,
\label{prop:tau2dg}
 $\tau_1(\delta) = [\min_{\gamma \in [0,\delta]} g(1,\gamma)]^{-1}$ and   $\tau_2(\delta) = [\min_{\gamma \in [0,\delta]} g(\gamma,1)]^{-1}$.
\end{prop}
\begin{proof}
 For $\gamma \in [0,1]$, the points $\left(1/g(1,\gamma), \gamma/g(1,\gamma)\right)$ lie on the curve $\{(x,y) \in [0,1]^2:g(x,y)=1, x\geq y\}$. The value $r_{1,\delta}$ is the maximum value of $1/g(1,\gamma)$ for $\gamma \in [0,\delta]$, hence $\tau_1(\delta) = [\min_{\gamma \in [0,\delta]} g(1,\gamma)]^{-1}$. A symmetric argument applies to $\tau_2(\delta)$.
\end{proof}
The right panel of Figure~\ref{fig:lambdatau} provides an illustration: in blue the value of $\delta$ is such that $\tau_1(\delta)<1$; in red the value of $\delta$ is such that $\tau_2(\delta)=1$. Further detail on this example is given in Section~\ref{sec:egg}.

The question arises: does Proposition~\ref{prop:tau2dg} still hold for $d>2, |C| = 1$? Let $g_C$ denote the gauge function for the limit set of $(X_{E,j}: j \in C)$. By Proposition~\ref{prop:marg}, we know that $g_{i,j}(x_i,x_j) = \min_{\bm{x}_{-i,-j} \in [0,\infty)^{d-2}}g_D(\bm{x})$. As such, equality will hold if $\arg\min_{\bm{x}_{-i,-j} \in [0,\infty)^{d-2}}g_D(\bm{x}) \in [0,\delta]^{d-2}$. Note that the dimension does indeed play a key role here: when looking at $\tau_j(\delta)$ for a $d$-dimensional problem, we are looking at the situation where $d-1$ coordinates are upper bounded by a growth rate determined by $\delta$. In contrast, when marginalizing and looking at $\tau_j(\delta)$ for a 2-dimensional problem, the $d-2$ coordinates that we have marginalized over are unrestricted and so can represent small or large values. As such, the answer to our question is negative in general.

Proposition~\ref{prop:taug} details the precise value of $\tau_C(\delta)$ in terms of $g$ for any dimension $d$. In a similar spirit to Section~\ref{sec:lambda}, define the boundary of the region $R_{C,\delta}$ as 
\[
B_{C,\delta} = B_{C,\delta}^1 \cup B_{C,\delta}^\delta,
\]
where 
\begin{align*}
B_{C,\delta}^1 &= \bigcup_{i \in C} \{\bm{x}\in\mathbb{R}^d_+: x_i  = 1, x_j \geq 1~ \forall j \in C \setminus i,~ x_k \leq \delta~ \forall k \in D \setminus C \}\\
B_{C,\delta}^\delta &= \bigcup_{i \in D \setminus C} \{\bm{x}\in\mathbb{R}^d_+: x_i  = \delta, x_j \leq \delta~ \forall j \in (D\setminus C) \setminus i,~ x_k \geq 1~ \forall k \in C \},
\end{align*}
so, for example, when $d=3$, 
\[B_{\{1,3\},\delta} = \{\bm{x}\in \mathbb{R}^3_+: x_1 = 1, x_2 \leq \delta, x_3 \geq 1\} \cup \{\bm{x}\in \mathbb{R}^3_+: x_1 \geq 1, x_2 = \delta, x_3 \geq 1\} \cup \{\bm{x}\in \mathbb{R}^3_+: x_1 \geq 1, x_2 \leq \delta, x_3 = 1\}.\]
 For $C=D$, $R_D = (1,\infty]^d$, and $B_D =\{\bm{x}: \min(\bm{x}) = 1\}$.
 
 \begin{prop}
 \label{prop:taug}
  	Assume the conditions of Proposition~\ref{prop:tau}. For any $C \subseteq D$, 
  \[\tau_C(\delta) = \left[\min_{\bm{y} \in B_{C,\delta}}g(\bm{y})\right]^{-1} =\left[\min_{\bm{y} \in B_{C,\delta}^1}g(\bm{y})\right]^{-1}. \]
 \end{prop}
\begin{proof}
 The vertex of the region $R_{C,\delta}$, or its boundary $B_{C,\delta}$, which has components 1 on the coordinates indexed by $C$, and $\delta$ in the other coordinates, lies on $\mathcal{S}_{\vee} :=\{\bm{x} \in \mathbb{R}^d_+: \max(\bm{x}) = 1\}$. The region $G \subseteq \mathcal{S}_{\vee}$, and because the coordinatewise supremum of $G$ is $\bm{1}$, the boundary of $G$ intersects with $\mathcal{S}_{\vee}$. Now consider scaling the region $R_{C,\delta}$ by $r_{C,\delta} \in (0,1]$ until it intersects with $G$. The point of intersection must lie on the boundary of the scaled region $r_{C,\delta}R_{C,\delta}$, i.e., on $r_{C,\delta}B_{C,\delta}$, and on the boundary of $G$, $\{\bm{x} \in \mathbb{R}^d_+: g(\bm{x}) = 1\}$. Therefore, there exists $\bm{x}^\star \in B_{C,\delta}$ such that $g(r_{C,\delta}\bm{x}^\star) = 1$, which is rearranged to give $\tau_C(\delta) = r_{C,\delta} = 1/ g(\bm{x}^\star)$. Furthermore, we must have that such a point $\bm{x}^\star = \arg\min_{\bm{y} \in B_{C,\delta}} g(\bm{y})$, otherwise there exists some $\bm{x}' \in B_{C,\delta}$ such that $g(\bm{x}')<g(\bm{x}^\star)$ and so $g(r_{C,\delta}\bm{x}') < 1$, meaning that $r_{C,\delta} \neq \min\{r \in (0,1]: r R_{C,\delta} \cap G = \emptyset\}$. We conclude that $\bm{x}^\star = \arg\min_{\bm{y} \in B_{C,\delta}} g(\bm{y})$, so $\tau_{C}(\delta) = 1/\min_{\bm{y} \in B_{C,\delta}} g(\bm{y})$.
 
 To show that $\arg\min_{\bm{y} \in B_{C,\delta}} g(\bm{y}) \in B_{C,\delta}^1$, let $\bar{\bm{x}} = \arg\min_{\bm{y} \in B_{C,\delta}^1} g(\bm{y})$, $\tilde{\bm{x}} = \arg\min_{\bm{y} \in B_{C,\delta}^\delta} g(\bm{y})$, and let $\tilde{x}_{l} = \min_{k \in C}\tilde{x}_k \geq 1$. Then $g(\tilde{\bm{x}}) = \tilde{x}_l g(\tilde{\bm{x}}/\tilde{x}_l)$, but $\tilde{\bm{x}}/\tilde{x}_l \in B_{C,\delta}^1$, so $g(\tilde{\bm{x}}/\tilde{x}_l) \geq g(\bar{\bm{x}})$ and hence $\tilde{x}_{l}g(\tilde{\bm{x}}/\tilde{x}_l) \geq g(\bar{\bm{x}})$  as $\tilde{x}_l \geq 1$.
\end{proof}

 When $d=2$ and $|C|=1$, we note that $B_{\{j\},\delta}^1 = \{\bm{x}:x_j = 1, x_i \leq \delta\}$, which gives the equality in Proposition~\ref{prop:tau2dg}.

\begin{proof}[Proof of Proposition~\ref{prop:lambdag}]
 The proof follows exactly as for the first equality in Proposition~\ref{prop:taug}, replacing $R_{C,\delta}, B_{C,\delta}$ and $r_{C,\delta}$ with $R_{\bm{\omega}}, B_{\bm{\omega}}$ and $r_{ \bm{\omega}}$.
\end{proof}

\subsubsection{Estimation of coefficients $\tau_C(\delta)$}
\label{sec:tauest}

When $C=D$, equation~\eqref{eq:SWT2} yields  $\pbb(\min_{i \in D} X_{P,i}> t) \in \RV^{\infty}_{-1/\tau_{D}}$, implying that $\tau_D$ can be estimated as the reciprocal of the tail index of the so-called structure variable $\min_{i \in D} X_{P,i}$. This is identical to estimating the residual tail dependence coefficient $\eta_D$, for which the Hill estimator is commonly employed. However, for $C$ with $|C|<d$, we assume $\pbb(\min_{i \in C} X_{P,i}> t , \max_{j \in D \setminus C} X_{P,j}< t^{\delta}) \in \RV^{\infty}_{-1/\tau_{C}(\delta)}$, but this representation does not lend itself immediately to an estimation strategy, as there is no longer a simple structure variable for which $1/\tau_{C}(\delta)$ is the tail index.

In order to allow estimation, \citet{Simpsonetal18} considered $\pbb(\min_{i \in C} X_{P,i}> t , \max_{j \in D \setminus C} X_{P,j}< (\min_{i \in C} X_{P,i})^{\delta})$, but they only offered empirical evidence that the assumed index of regular variation for this probability was the same as in equation~\eqref{eq:SWT2}. We now prove this to be the case.

Define $R_{C,\delta}^x = \{\bm{x}\in \mathbb{R}^d_+: x_i>1, i \in C,~x_j \leq \delta \min_{l \in C} x_l, j \in D\setminus C\}$, and $B_{C,\delta}^x = B_{C,\delta}^{x,1} \cup B_{C,\delta}^{x,\delta}$ to be its boundary, where
\begin{align*}
 B_{C,\delta}^{x,1}= B_{C,\delta}^1 &= \bigcup_{i \in C} \{\bm{x}\in\mathbb{R}^d_+: x_i  = 1, x_j \geq 1~ \forall j \in C \setminus i,~ x_k \leq \delta \min_{l \in C} x_l = \delta~ \forall k \in D \setminus C \} \\
 B_{C,\delta}^{x,\delta} &=  \bigcup_{i \in D \setminus C} \{\bm{x}\in\mathbb{R}^d_+: x_i  = \delta\min_{l \in C} x_l,~ x_j \leq \delta\min_{l \in C} x_l~ \forall j \in (D\setminus C) \setminus i,~ x_k \geq 1~ \forall k \in C \},
\end{align*}
and we specifically note the equality $ B_{C,\delta}^{x,1} =  B_{C,\delta}^{1}$. 

\begin{prop} Assume the conditions of Proposition~\ref{prop:tau}. If $\pbb(\min_{i \in C} X_{P,i}> t , \max_{j \in D \setminus C} X_{P,j}< t^{\delta}) \in \RV^{\infty}_{-1/\tau_{C}(\delta)}$, and $\pbb(\min_{i \in C} X_{P,i}> t , \max_{j \in D \setminus C} X_{P,j}< (\min_{i \in C} X_{P,i})^{\delta}) \in \RV_{-1/\tilde{\tau}_C(\delta)}^\infty,$ then $\tilde{\tau}_C(\delta) = \tau_C(\delta)$.
\end{prop}
\begin{proof}
Define $r_{C,\delta}^x = \min\left\{r \in [0,1] : r R_{C,\delta}^x \cap G =\emptyset \right\}$, where for $r>0$, $rR_{C,\delta}^x = \{\bm{x}\in \mathbb{R}^d_+: x_i>r, i \in C,~x_j \leq \delta \min_{l \in C} x_l, j \in D\setminus C\}$. Similarly to Propositions~\ref{prop:lambda} and~\ref{prop:tau},
  \begin{align*}
-\log \pbb(\bXe \in tr_{C,\delta}^x R_{C,\delta}^x) = -\log \pbb(\min_{i \in C} X_{E,i}>tr_{C,\delta}^x, \max_{j \in D \setminus C} X_{E,j}< \delta \min_{i \in C} X_{E,i}) \sim t,
 \end{align*}
 and we conclude $\tilde{\tau}_C(\delta) = r_{C,\delta}^{x}$. As in Proposition~\ref{prop:taug}, we have $r_{C,\delta}^{x}=\min_{\bm{y} \in B_{C,\delta}^x} g(\bm{y})$.  Noting again that $\arg\min_{\bm{y} \in B_{C,\delta}^x} g(\bm{y}) \in B_{C,\delta}^{x,1}=B_{C,\delta}^1$ shows that $r_{C,\delta}^{x} = r_{C,\delta} = \tau_{C}(\delta)$.
\end{proof}

\section{Examples}
\label{sec:examples}

We illustrate several of the findings of Section~\ref{sec:connections} with some concrete examples. In Section~\ref{sec:egd2} we focus on the intuitive and geometrically simple case $d=2$; in Section~\ref{sec:egd3}, we examine some three-dimensional examples for which visualization is still possible but more intricate. \iftoggle{arxiv}{}{Additional examples are given in the arXiv version of this article.}

Proposition~\ref{prop:logdens} implies that on $\mathbb{R}^d_+$, the same limit set $G$ as in exponential margins will arise for any marginal choice with $\psi_j(x) \sim x$, $x \to \infty$, provided $e^{-\psi_j(x)}$ is a von Mises function. In some of the examples below, it is convenient to establish a limit set and its gauge function using this observation rather than transforming to exactly exponential margins.

Models with convenient dependence properties are often constructed through judicious combinations of random vectors with known dependence structures; see, for example, \citet{Engelkeetal18} for a detailed study of so-called random scale or random location constructions. In Section~\ref{sec:mixing}, we use our results to elucidate the shape of the limit set when independent exponential-tailed variables are mixed additively. The spatial dependence model of \citet{HuserWadsworth18} provides a case study.

\subsection{Examples and illustrations for $d=2$}
\label{sec:egd2}
All of the examples considered in this section are symmetric, so, for the conditional extremes representation and coefficients $\tau_j(\delta)$, we only consider one case, omitting the subscript on the quantities~$\alpha_j$ and~$\beta_j$. Table~\ref{tab:bivariate} summarizes the dependence information from various bivariate distributions described in \iftoggle{arxiv}{Sections~\ref{exGpos}--\ref{se:eghr}.}{Sections~\ref{exGpos}--\ref{sec:egievd}, as well as the arXiv version.}

\afterpage{%
	\clearpage
	\begin{landscape}
		\centering 
\begin{tabular}[htbp]{llllll}
	\hline
	Copula & $g(x,y)$ & $\lambda(\omega,1-\omega)$ & $\eta_{1,2}$ & $\tau_{1}(\delta)$ & $\alpha,\beta$\\
	\hline
	
	Gaussian $\rho\geq 0$ & $\frac{x + y -2\rho(xy)^{1/2}}{1-\rho^2}$ &  $\begin{cases}
	\max(\omega,1-\omega), & \frac{\min(\omega,1-\omega)}{\max(\omega,1-\omega)} \leq \rho^2\\
	\frac{1 -2\rho(\omega(1-\omega))^{1/2}}{1-\rho^2}, &  \frac{\min(\omega,1-\omega)}{\max(\omega,1-\omega)} \geq \rho^2
	\end{cases}$ & $\frac{1+\rho}{2}$ & $\begin{cases}1 & \delta \geq \rho^2\\ \frac{1-\rho^2}{1+\delta-2\rho\delta^{1/2}} & \delta < \rho^2 \end{cases}$ & $\begin{array}{l} \alpha=\rho^2 \\ \beta=1/2\end{array}$ \\\hline
	
	Logistic GP & $\frac{1}{\theta}\max(x,y) + (1-\frac{1}{\theta})\min(x,y)$ & $\max(\omega,1-\omega)$ & 1 & $[\theta^{-1}+(1-\theta^{-1}\delta)]^{-1}$ & $\begin{array}{l} \alpha=1 \\ \beta=0\end{array}$\\\hline
	
		\shortstack{Inverted\\ logistic}& $(x^{1/\theta}+y^{1/\theta})^\theta $ & $g(\omega,1-\omega)$ & $2^{-\theta}$ & $1$ & $\begin{array}{l} \alpha=0 \\ \beta=1-\theta\end{array}$\\\hline

	\shortstack{H\"{u}sler--Reiss\\ GP} &  $\begin{cases}
	\infty, & x \neq y\\
	x, &  x=y
	\end{cases}$ & $\max(\omega,1-\omega)$ & 1 & $\begin{cases}
	0, & \delta<1\\
	1, &  \delta=1
	\end{cases}$ & $\begin{array}{l} \alpha=1 \\ \beta \mbox{ undetermined}\end{array}$\\\hline		
		
	\shortstack{Inverted\\ H\"{u}sler--Reiss}& $ x\Phi\left(\frac{\lambda}{2} + \frac{\log\left(\frac{x}{y}\right)}{\lambda}\right) + y\Phi\left(\frac{\lambda}{2} + \frac{\log\left(\frac{y}{x}\right)}{\lambda}\right) $ & $g(\omega,1-\omega)$ & $\left[2\Phi\left(\frac{\lambda}{2}\right)\right]^{-1} $ & $1$ & $\begin{array}{l} \alpha=0 \\ \beta=1\end{array}$\\\hline
		
\end{tabular}
		\captionof{table}{Summary of dependence measures across a range of bivariate examples}
		\label{tab:bivariate}
	\end{landscape}
	\clearpage
}

\subsubsection{Meta-Gaussian distribution: nonnegative correlation}\label{exGpos}

Starting with a Gaussian bivariate random vector and transforming its margins to standard exponential, we obtain a meta-Gaussian distribution with exponential margins. Such a distribution inherits the copula of the Gaussian distribution. For simplicity, we consider the case where the underlying Gaussian random vector has standard normal components with correlation $\r$.  

Then, for $\r\ge0$, the joint probability density $f_E$ satisfies:
\begin{align*}
    -\log f_E(tx,ty)/t = (x + y -2\rho(xy)^{1/2})/(1-\rho^2) + O(\log t/t),\quad t\to\nf,\quad x,y\ge 0, 
\end{align*}
so that $g(x,y) = (x + y -2\rho(xy)^{1/2})/(1-\rho^2)$. The convergence in \eqref{q1} holds on $[0,\nf)^d$ and hence the limit set exists and is given by $\{\xb\in[0,\nf)^d: g(\xb)\le1\}$. This is example~(ii) in Figure~\ref{fig:2dgauge}.

\paragraph{Conditional extremes:} Setting $g(\alpha,1) = 1$ leads to $(\alpha^{1/2} - \rho)^2 = 0$, i.e., $\alpha = \rho^2$. For $\beta$ we have $g(\rho^2+u, 1) - 1 = u^2/\{2\rho(1-\rho^2)\} + O(u^3) \in \RV_2^0$, hence $\beta = 1/2$.

\paragraph{Function $\lambda(\bm{\omega})$:} By Proposition~\ref{prop:lambdag}, we need to find $1/r_{\omega} = \min_{\bm{y} \in B_{\omega}} g(x,y)$. If $\min(\omega,1-\omega)/\max(\omega,1-\omega) \leq \rho^2$, then $\min_{\bm{y} \in B_{\omega}} g(x,y) = 1$, with the minima occuring at the points $(1,\rho^2)$, $(\rho^2,1)$. Otherwise, if $\min(\omega,1-\omega)/\max(\omega,1-\omega) \geq \rho^2$, then $\min_{\bm{y} \in B_{\omega}} g(x,y) = g(1, \min(\omega,1-\omega)/\max(\omega,1-\omega))$. Putting this together with Proposition~\ref{prop:lambda}, we find
\[
 \lambda(\omega,1-\omega) = \begin{cases}
                    \max(\omega,1-\omega), & \min(\omega,1-\omega)/\max(\omega,1-\omega) \leq \rho^2\\
                    g(\omega,1-\omega) = \frac{1 -2\rho(\omega(1-\omega))^{1/2}}{1-\rho^2}, &  \min(\omega,1-\omega)/\max(\omega,1-\omega) \geq \rho^2.
                   \end{cases}
\]
This is the same form as given in \citet{WadsworthTawn13}. We therefore have $\eta_{1,2} = [2g(1/2,1/2)]^{-1} = g(1,1)^{-1} = (1+\rho)/2$.

\paragraph{Coefficients $\tau_j(\delta)$:} From Proposition~\ref{prop:tau2dg}, we have $\tau_1(\delta) = [\min_{\gamma \in [0,\delta]} g(1,\gamma)]^{-1} = [g(1,\min(\delta,\rho^2))]^{-1}$. Therefore, $\tau_{1}(\delta) = 1$ if $\delta \geq \rho^2$, else $\tau_{1}(\delta) = (1-\rho^2)/(1+\delta-2\rho\delta^{1/2}) < 1$. Note that these values are very laborious to calculate via Gaussian survival functions, and they were not given in \citet{Simpsonetal18}.

\subsubsection{Meta-Gaussian distribution: negative correlation}\label{exGneg}
When $\rho<0$, Proposition~\ref{prop:logdens} cannot be applied as the continuous convergence condition~\eqref{q1} does not hold along the axes. Hence, we only gain a partial specification, when $x>0, y>0$, through this route. Instead, here we can apply Proposition~\ref{prop:DMR} since the limit function $g$ in~\eqref{eq:DMR} satisfies the monotonicity condition given immediately thereafter. This limit function is given by
\begin{align*}
 g(x,y) = \begin{cases}
           (x + y -2\rho(xy)^{1/2})/(1-\rho^2), & x>0,~ y>0,\\
           x, & y=0,\\
           y, & x=0.
          \end{cases}
\end{align*}

\begin{figure}[]
 \includegraphics[width=1\textwidth]{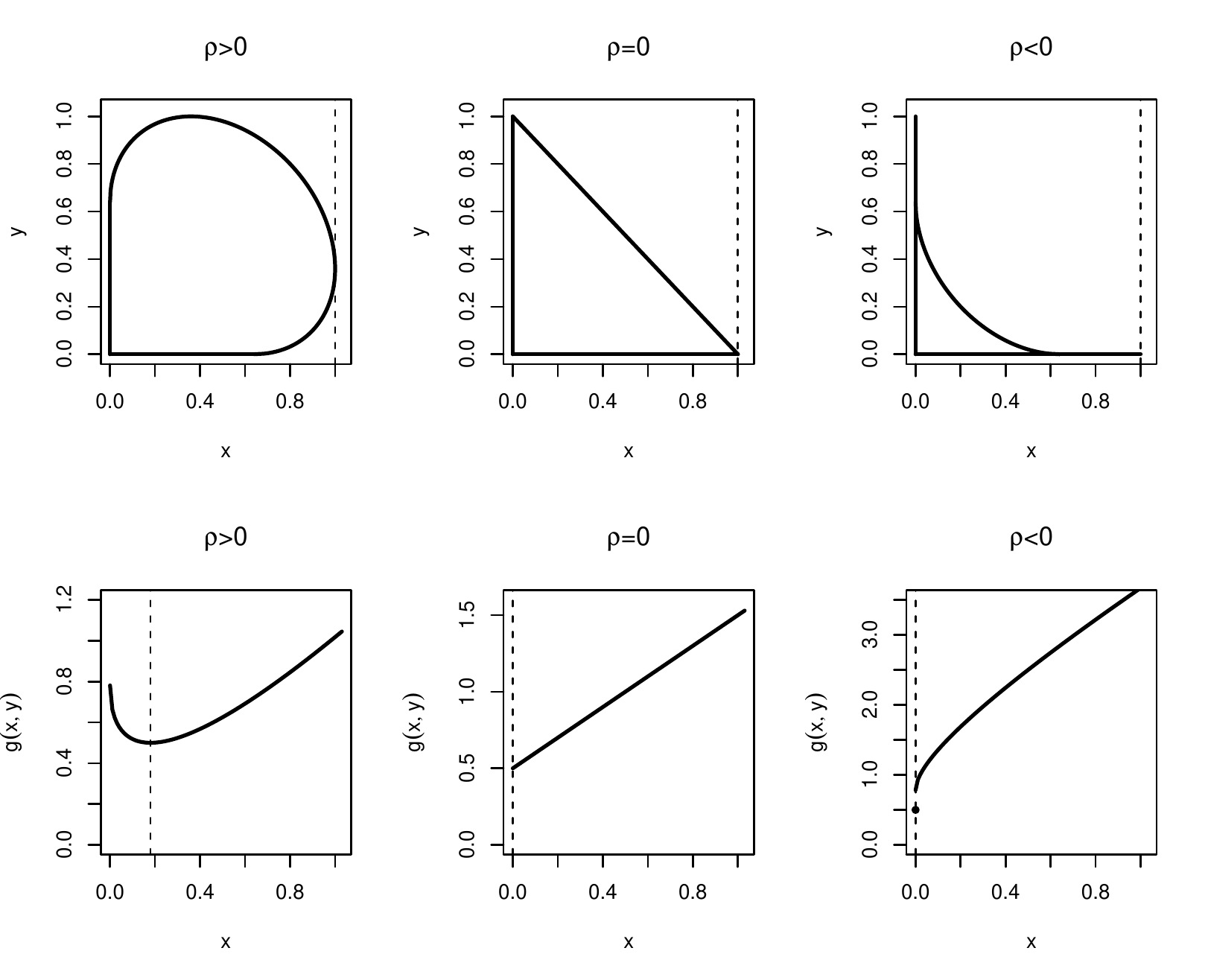}
 \caption{Top row: Limit sets $\{(x,y): g(x,y)\le1\}$ for a bivariate meta-Gaussian distribution with exponential margins. Bottom row: corresponding plots of function $g(x,y)$ for a fixed value of $y$.}
 \label{fig:2dgauss}
\end{figure}

Figure~\ref{fig:2dgauss} illustrates the limit sets $G$ for the three cases $\rho>0,\rho=0$ and $\rho<0$. In the latter case, large values of one variable tend to occur with small values of the other, which causes the limit set to include lines along the axes and the function $g$ is not continuous. Such difficulties can be alleviated by consideration of Laplace margins for distributions displaying negative dependence, which is discussed further in Section~\ref{sec:Discussion}.

\subsubsection{Logistic generalized Pareto copula}
\label{sec:eglog}

The logistic generalized Pareto distribution with conditionally exponential margins ($\pbb(X_{\tilde{E}}>x) = \pbb(X_{\tilde{E}}>0)e^{-x}$, $x>0$) and dependence parameter $\theta \in(0,1)$ satisfies
\begin{align*}
   f_{\tilde{E}}(x,y) &= \theta^{-1}2^{-\theta} e^{-(x+y)/\theta}(e^{-x/\theta}+e^{-y/\theta})^{\theta-2}\\
   -\log f_{\tilde{E}}(tx,ty)/t &=  \theta^{-1}\max(x,y) + (1-\theta^{-1})\min(x,y) + O(1/t),
\end{align*}
so the gauge function is $g(x,y) = \theta^{-1}\max(x,y) + (1-\theta^{-1})\min(x,y)$. This form of gauge function is found throughout several symmetric asymptotically dependent examples, such as those distributions whose spectral measure $H$ places no mass on 0 and 1 and possess densities that are regularly varying at the endpoints 0, 1 such that $\mathrm{d}H(w)/\mathrm{d}w \in \RV_{1/\theta-2}^{0}$, $-\mathrm{d}H(1-w)/\mathrm{d}w \in \RV_{1/\theta-2}^{0}$. This is example~(i) in Figure~\ref{fig:2dgauge}.

\paragraph{Conditional extremes:} Solving for $g(\alpha,1)=1$, we obtain $\alpha=1$, whilst $g(1+u,1) -1 = u/\theta \in \RV_1^0$, hence $\beta=0$.

\paragraph{Function $\lambda(\bm{\omega})$:} We have that $\arg\min_{\bm{y} \in B_\omega}g(\bm{y}) = (1,1)$, so $r_{ \omega} = 1$ and $\lambda(\omega,1-\omega) = \max(\omega,1-\omega)$. Therefore $\eta_{1,2} =1$.

\paragraph{Coefficients $\tau_j(\delta)$:}  $\tau_1(\delta) = [\min_{\gamma \in [0,\delta]} g(1,\gamma)]^{-1} = [g(1,\delta)]^{-1} = [\theta^{-1}+(1-\theta^{-1}\delta)]^{-1}$. This matches the value calculated in the Supplementary Material of \citet{Simpsonetal18}.

\subsubsection{Inverted extreme value distribution}
\label{sec:egievd}
The inverted extreme value copula is the joint lower tail of an extreme value copula, translated to be the joint upper tail. That is, if $(U_1,U_2)$ have an extreme value copula with uniform margins, then $(1-U_1,1-U_2)$ have an inverted extreme value copula. In two dimensions, its density in exponential margins may be expressed as
\[
f_E(x,y) = \{l_{1}(x,y)l_2(x,y)-l_{12}(x,y)\}\exp\{-l(x,y)\},
\]
where $l(\bm{x}) = V(1/\bm{x})$, for $V$ the exponent function in~\eqref{eq:V}, is the 1-homogeneous stable tail dependence function \citep[e.g.,][Ch.8]{Beirlantetal04} of the corresponding extreme value distribution, and $l_1(x,y) = \partial l(x,y) /\partial x$, etc. We thus have
\[
-\log f_E(tx,ty)/t = l(x,y) + O(1/t),\qquad t \to\infty, 
\]
so $g(x,y) = l(x,y)$. 

\paragraph{Conditional extremes:} Stable tail dependence functions always satisfy $l(x,0) = x, l(0,y) = y$ and so $g(1,0) = g(0,1) = 1$. Hence, if $\alpha=0$ is the only solution to $g(\alpha,1) =1$, then $a(x)/x \sim 0$. An example of this is given by the inverted extreme value logistic copula, whereby $l(x,y) = (x^{1/\theta}+y^{1/\theta})^\theta, \theta \in (0,1]$. This is example~(iii) of Figure~\ref{fig:2dgauge}, for which we have $\alpha=0$ and $\beta=1-\theta$.

\iftoggle{arxiv}{\citet{PapastathopoulosTawn16} study conditional extreme value limits for general classes of inverted extreme value distributions. One case they consider is where the spectral measure defining the distribution in two dimensions has support on a sub-interval $[w_l,w_u] \subset [0,1]$. We give a simple example to illustrate their findings in this context. Let $w_l=1-w_u = w^\star \in [0,1/2)$, and $H(w) = (w-w^\star)/(1-2w^\star), w \in [w^\star,1-w^\star]$. Then
\begin{align*}
    l(x,y) &= \int_{w^\star}^{1-w^\star} \max(wx,(1-w)y)\mathrm{d}H(w) \\&= \begin{cases}
    x, & y< x w^\star/(1-w^\star)\\
    \frac{x}{1-2w^\star}\left\{(1-w^\star)^2-\left(\frac{y}{x+y}\right)^2\right\}+\frac{y}{1-2w^\star}\left\{(1-w^\star)^2-\left(\frac{x}{x+y}\right)^2\right\}, &  x w^\star/(1-w^\star) \leq y \leq x (1-w^\star)/w^\star\\
    y,& y> x (1-w^\star)/w^\star.
    \end{cases}
\end{align*}
In this case, $g(\alpha,1) = 1$ for all $\alpha \in [0,w^\star/(1-w^\star)]$, so by Proposition~\ref{prop:CEgauge}~(iii), $\alpha=w^\star/(1-w^\star)$. Following a Taylor expansion in which linear terms in $u$ vanish, we have $g(w^\star/(1-w^\star) + u, 1) - 1 \in \RV^0_{2}$, so $\beta=1/2$. \citet{PapastathopoulosTawn16} show that for this distribution, $a(x) = x w^\star/(1-w^\star)$ and $b(x)=x^{1/2}$.

A further interesting example studied by \citet{PapastathopoulosTawn16} is that of the inverted H\"{u}sler--Reiss distribution, for which
\[
g(x,y) = x\Phi(\lambda/2 + \log(x/y)/\lambda) + y\Phi(\lambda/2 + \log(y/x)/\lambda), \qquad \lambda>0.
\]
The unique solution to $g(\alpha,1) = 1$ is $\alpha=0$ (obtained as a limit) and $g(u,1)-1$ is rapidly varying, i.e., $g(u,1)-1 \in \RV_{\infty}^0$, corresponding to $\beta=1$. The forms of the normalization functions detailed in \citet{PapastathopoulosTawn16} are
\begin{align*}
a(x) = x \exp\left\{-\lambda(2\log x)^{1/2} +\frac{\lambda\log\log x}{(2\log x)^{1/2}}+\lambda^2/2\right\},\qquad b(x) = a(x)/(\log x)^{1/2},
\end{align*}
for which $a(x)/x \to 0$ and $b(x) \in \RV_1^{\infty}$.}{Further examples in this class are given in the arXiv version.}

\paragraph{Function $\lambda(\bm{\omega})$:} Since $g(x,y) = l(x,y)$, and $l$ is a convex function satisfying $l(x,0)=x,l(0,y) = y$, $\arg\min_{\bm{y} \in B_{\omega}} g(\bm{y}) = (\omega,1-\omega)/\max(\omega,1-\omega)$. Hence, $\lambda(\omega,1-\omega) = g(\omega,1-\omega)$ in this case. 

\paragraph{Coefficients $\tau_j(\delta)$:} Since $g(1,0)=1$, we have $\tau_1(\delta) = 1$ for all $\delta \in [0,1]$.

\iftoggle{arxiv}{
\subsubsection{H\"{u}sler--Reiss generalized Pareto copula}
\label{se:eghr}
The bivariate H\"{u}sler--Reiss generalized Pareto distribution with conditionally exponential margins has density
\[
f_{\tilde{E}}(x,y) ~~\propto~~ \exp\left\{-\frac{1}{2}\left(\frac{(x-y)^2}{2(1-\rho)}+x+y\right)\right\},
\]
from which it can be seen that
\[
-\log f_{\tilde{E}}(t\bm{x})/t \to \begin{cases}
\infty, & x \neq y,\\
x, & x=y.
\end{cases}
\]
While Proposition 2 cannot be applied here due to the lack of uniform convergence, the form of the limit set is nonetheless $G=\{(x,y):x=y \leq 1\}$, which is the same limit set as arises under perfect dependence. This can be explained by the construction of the H\"{u}sler--Reiss model, which has a dependence structure asymptotically equivalent to that of $E+(Z_1,Z_2)$, where $E \sim \mbox{Exp}(1)$ is independent of $(Z_1,Z_2) \sim N_2(\bm{0},\Sigma)$. The marginal distributions of such a construction satisfy $\psi(x) \sim x$, $x \to \infty$ so that with the sample cloud constructed by taking $r_n = \log n$, the contribution from the Gaussian component converges in probability to zero, leaving only the contribution from the common term $E$.

\paragraph{Conditional extremes:} We have $g(\alpha,1)=1$ for $\alpha=1$, but cannot use Proposition~\ref{prop:CEgauge} to determine $\beta$ since $g(1+u,1) -1 \not\in \RV^0$.

\paragraph{Function $\lambda(\bm{\omega})$:} The quantity $r_{\omega}=1$ and so $\lambda(\omega,1-\omega)=\max(\omega,1-\omega)$, and $\eta_{1,2}=1$.

\paragraph{Coefficients $\tau_j(\delta)$:} $\tau_1(\delta)=[\min_{\gamma \in [0,\delta]} g(1,\gamma)]^{-1} = 0$ for $\delta<1$ and $=1$ for $\delta=1$. This implies that $\pbb(X_P>t,Y_P \leq t^{\delta}) \in \RV^\infty_{-\infty}$, i.e., is rapidly varying.
}{}

\subsubsection{Density defined by $g$}
\label{sec:egg}
If $g:\mathbb{R}^d_+ \to \mathbb{R}_+$ is a gauge function describing a limit set $G$, then $f(\bm{x}) = e^{-g(\bm{x})}/(d! |G|)$ is a density \citep[see][]{BalkemaNolde10}. In general, except for the case of $g(\bm{x}) = \sum_{i=1}^d x_i$, the margins are not exactly exponential, and may be heavier than exponential, for example in the case $g(\bm{x}) = \max_{1\leq i \leq d}(x_i)$. 

We consider the density defined by $g(x,y) = \max\{(x-y)/\theta,(y-x)/\theta,(x+y)/(2-\theta)\}$, $\theta \in (0,1]$: this is example~(vi) in Figure~\ref{fig:2dgauge}, and illustrated in Figure~\ref{fig:lambdatau}. The marginal density is given by
\[
 [2 e^{-x} - \theta e^{-x/\theta} -2(1-\theta)e^{-x/(1-\theta)}]/[4\theta - 3\theta^2].
 \]

\paragraph{Conditional extremes:} Solving for $g(\alpha,1)=1$, we obtain $\alpha=1-\theta$, whilst $g(1-\theta+u,1) -1 = u/(2-\theta) \in \RV_1^0$, hence $\beta=0$.

\paragraph{Function $\lambda(\bm{\omega})$:} If $\min(\omega,1-\omega)/\max(\omega,1-\omega) \leq 1-\theta$, then $\arg\min_{\bm{y} \in B_\omega} = (1,1-\theta)$, or $(1-\theta,1)$ and $r_{\omega}=1$; otherwise, $\arg\min_{\bm{y} \in B_\omega} = (1,\omega/(1-\omega))$ or $((1-\omega)/\omega, 1)$, and $r_{\omega} =\{1+ \min(\omega,1-\omega)/\max(\omega,1-\omega)/(2-\theta)\}$. As such
\[
 \lambda(\omega,1-\omega) = \begin{cases}
                    \max(\omega,1-\omega), & \min(\omega,1-\omega)/\max(\omega,1-\omega) \leq 1-\theta\\
                    g(\omega,1-\omega) = \frac{1}{2-\theta}, &  \min(\omega,1-\omega)/\max(\omega,1-\omega) \geq 1-\theta,
                   \end{cases}
\]
and the residual tail dependence coefficient $\eta_{1,2} = 1-\theta/2$.

\paragraph{Coefficients $\tau_j(\delta)$:} $\tau_1(\delta) = [\min_{\gamma \in [0,\delta]} g(1,\gamma)]^{-1} = [g(1,\min(\delta,1-\theta))]^{-1}$. Therefore, $\tau_{1}(\delta) = 1$ if $\delta \geq 1-\theta$, else $\tau_{1}(\delta) = \theta/(1-\delta)< 1$.

\iftoggle{arxiv}{
\subsubsection{Boundary case between asymptotic dependence and independence}
We give two examples of distributions whose limit set is described by the gauge function $g(x,y)=\max(x,y)$. The first of these displays asymptotic dependence, i.e., $\nu(\mathbb{E}_{1,2})>0$, for $\nu$ as in Section~\ref{sec:MRVintro}, while the second displays asymptotic independence. In both cases, $\lambda(\bm{\omega})=\max(\omega,1-\omega), \eta_{1,2}=1$, and $\tau_{1}(\delta)=\tau_2(\delta)=1$ for all $\delta \in (0,1)$. The conditional extremes convergences are however rather different, and need to be derived carefully as some of the hypotheses in Proposition~\ref{prop:CEgauge} fail.

\paragraph{Example 1: asymptotic dependence}

We consider a particular instance of a bivariate extreme value distribution with spectral measure density $\mathrm{d}H(w)/\mathrm{d}w = h(w) = C w^{-1}(1-w)^{-1}\exp\{-(-\log w)^{1/2}-(-\log(1-w))^{1/2}\}$. The density is regularly varying at 0 and 1 with index $-1$, but its integral is finite. For the exponent function $V$ as in~\eqref{eq:V}, $h(w)=-V_{12}(w,1-w)$ with $V_{12}(x,y)=\partial^2 V(x,y)/(\partial x\partial y)$, so 
\[
h\left(\frac{x}{x+y}\right) = -V_{12}\left(\frac{x}{x+y}, \frac{y}{x+y}\right) = -(x+y)^3 V_{12}(x,y), 
\]
and therefore
\[
-V_{12}(x,y) = \frac{C}{(x+y)xy}\exp\{-(-\log x + \log(x+y))^{1/2}-(-\log y + \log(x+y))^{1/2}\}.
\]
For the corresponding multivariate max-stable or generalized Pareto distribution with Gumbel/exponential margins, the gauge function is determined by 
\[
f_{\tilde{E}}(x,y) ~~\propto~~ e^{x+y}V_{12}(e^x,e^y) ~~\propto~~ \frac{1}{e^x+e^y} \exp\{-(- x + \log(e^x+e^y))^{1/2}-(-y + \log(e^x+e^y))^{1/2}\},
\]
and $-\log f_{\tilde{E}}(tx,ty) \sim t\max(x,y)$. That is, $g(x,y) = \max(x,y)$.

We have $g(\alpha,1)=1$ for all $\alpha \in [0,1]$, so $\alpha=1$ is the maximum such value, and $g(1,1+u)-1=u \in \RV^0_1$, giving $\beta=0$. Taking $a^j(t) = t$, $b^j(t) = 1$, we see
\[
f_{\tilde{E}}(t,t+z)e^t ~~\propto~~ (1+e^z)^{-1} \exp\{-[\log(1+e^{-z})]^{1/2}\}, \qquad z \in \mathbb{R}.
\]
The limit is in fact exact for all $t$ here because of the multivariate generalized Pareto form.

\paragraph{Example 2: asymptotic independence}

Consider the density $f(x,y) = e^{-\max(x,y)}/2$, $x,y>0$. The marginal densities are $f(x)=(1+x)e^{-x}/2$, which is a mixture of $\rm{Gamma(1,1)}$ and $\rm{Gamma(2,1)}$ densities, and is heavier tailed than standard exponential. We firstly verify asymptotic independence by examining $\lim_{t \to \infty}\pbb(X>t|Y>t)$. The joint survival function $\pbb(X>x,Y>y)=\{2+\max(x,y)-\min(x,y)\}e^{-\max(x,y)}/2$, while $\pbb(X>x)=e^{-x}(1+x/2)$. Hence $\lim_{t \to \infty} \pbb(X>t|Y>t) = 0$.

Again, Proposition~\ref{prop:logdens} ensures that the limit set with gauge function $g(x,y)=\max(x,y)$ would also arise in exponential margins. However, we make the transformation explicitly here to study the conditional extremes convergence on the same scale as previously. The change to exponential margins entails $X_{E} = X-\log(1+X/2)$, leading to the density
\[
f_E(x,y) = \exp[-\max\{x + \log (1+x/2+O(\log x)), y + \log (1+y/2+O(\log y))\}][2^{-1}+O(x^{-1})+O(y^{-1})], 
\]
from which we can also see $-\log f_E(tx,ty) \sim t \max(x,y)$.

If we were to suppose that a conditional extremes limit with support including $(0,\infty)$ exists, then Proposition~\ref{prop:CEgauge}~(i) and~(iii) would mean $\alpha_1=\alpha_2=1$, while $g(1+u,1)-1 =u \in \RV^0_1$, so that Proposition~\ref{prop:CEgauge}~(ii) would give $\beta_1=\beta_2=0$. However, for positive $z$, consideration of $b(t)f_E(t,a(t)+b(t)z)e^t$ yields no possibilities for a non-degenerate limit. Nonetheless, for $z<0$ we can take $a(t)=b(t)=t$ leading to 
\[
tf_E(t,t+tz)e^t \sim (1+t/2)^{-1}t/2 \to 1, \qquad t \to \infty, \qquad z \in (-1,0),
\]
i.e., a uniform limit distribution on $(-1,0)$.

We comment that the results of Proposition~\ref{prop:CEgauge} do indeed focus predominantly on the positive end of the support for limit distributions, but most known examples of conditional limits have support including $(0,\infty)$.  A natural next step is to consider the implications relating to negative support. We particularly note the possibility that the order of regular variation of the two functions $g(1,\alpha_1+u)-1 \in \RV^0_{1/(1-\beta_1^+)}$ and $g(1,\alpha_1-u)-1\in \RV^0_{1/(1-\beta_1^-)}$ need not be equal, though for each of our examples where both functions are regularly varying, $\beta^+=\beta^-$. If $\beta^+>\beta^-$, it seems likely that a limit distribution with positive support only would arise, and vice versa when $\beta^+<\beta^-$.}{Further examples presented in the arXiv version include the H\"{u}sler--Reiss generalized Pareto copula, and a two boundary cases with $g(x,y)=\max(x,y)$, one displaying asymptotic dependence and one asymptotic independence. In the latter of these, we find that there is no conditional extreme value limit with positive support, but there is one with negative support. We comment that the results of Proposition~\ref{prop:CEgauge} do indeed focus predominantly on the positive end of the support for limit distributions, but most known examples of conditional limits have support including $(0,\infty)$.  A natural next step is to consider the implications relating to negative support. We particularly note the possibility that the order of regular variation of the two functions $g(1,\alpha_1+u)-1 \in \RV^0_{1/(1-\beta_1^+)}$ and $g(1,\alpha_1-u)-1\in \RV^0_{1/(1-\beta_1^-)}$ need not be equal, though for each of our examples where both functions are regularly varying, $\beta^+=\beta^-$. If $\beta^+>\beta^-$, it seems likely that a limit distribution with positive support only would arise, and vice versa when $\beta^+<\beta^-$.}

\subsection{Examples and illustrations for $d=3$}
\label{sec:egd3}
In this section we give two examples, focusing on issues that arise for $d>2$.
\subsubsection{Gaussian copula}
\begin{figure}[]
 \includegraphics[width=0.33\textwidth]{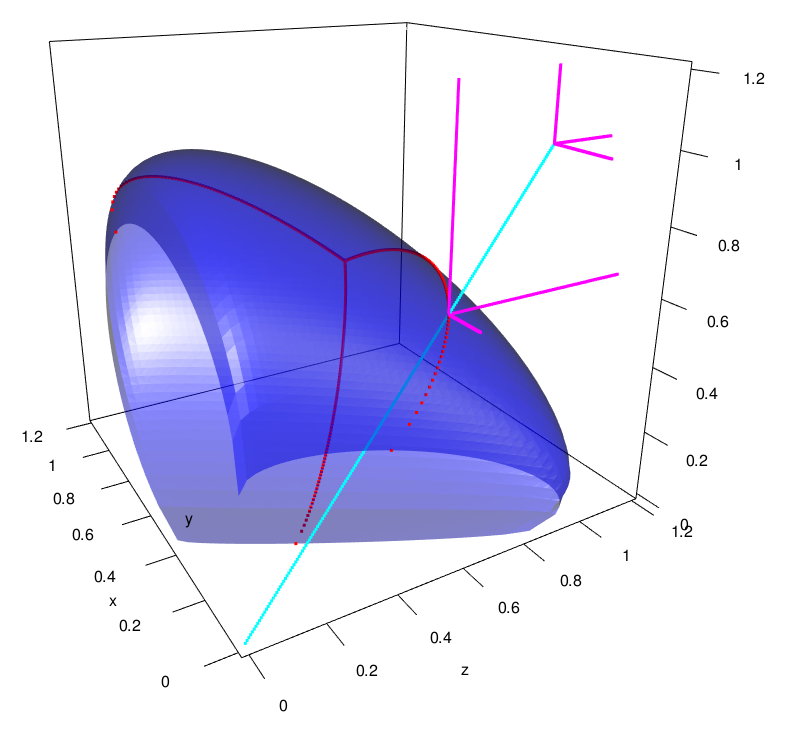}
 \includegraphics[width=0.33\textwidth]{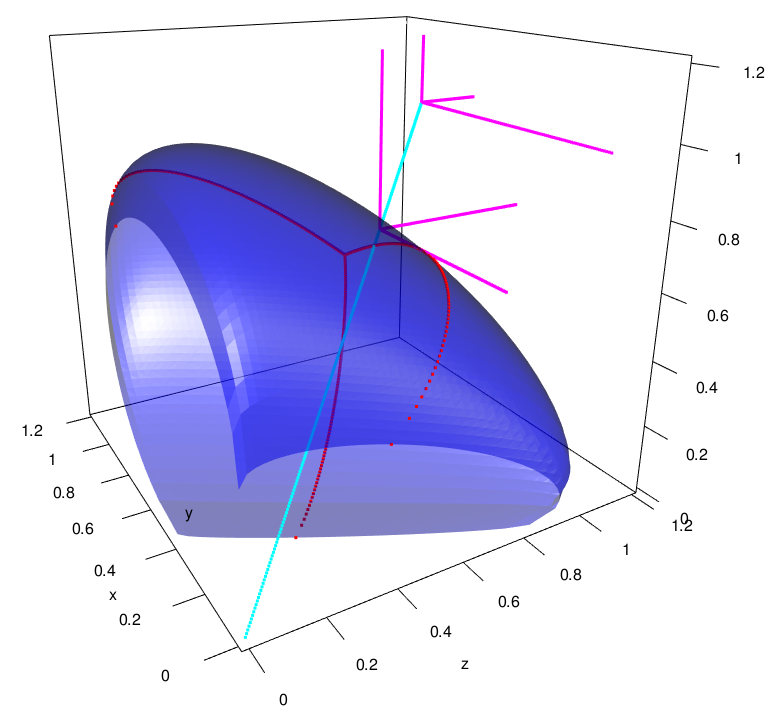}
  \includegraphics[width=0.33\textwidth]{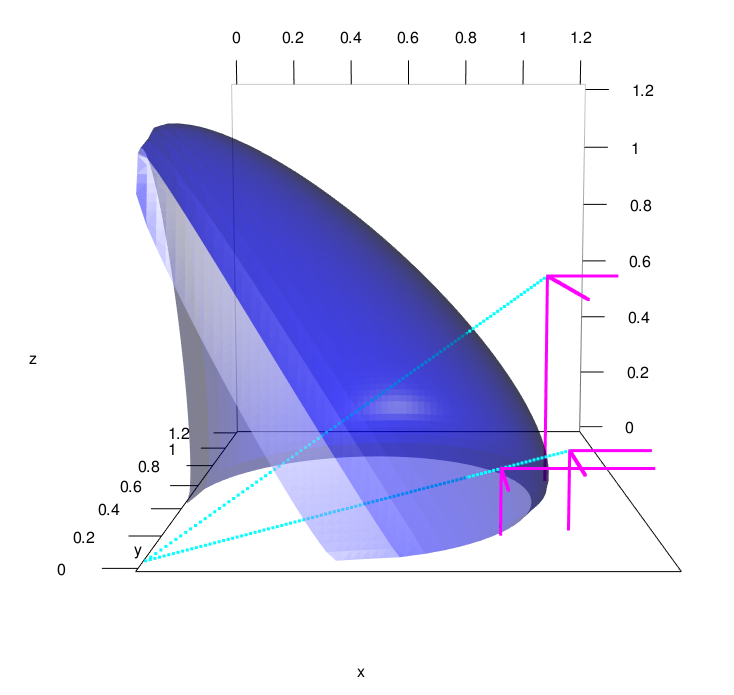}
\caption{Level set $g(\bm{x})=1$ for a trivariate meta-Gaussian distribution with exponential margins. The left panel illustrates $\tau_{2,3}(0.2)$: the boundary set indicated by the pink lines is scaled along the blue trajectory until it touches $G$, which happens in this case at the corner point $(0.2,1,1)/g(0.2,1,1)$. The centre panel illustrates $\tau_{2,3}(0.8)$: the boundary set is again pulled back along the indicated trajectory until it touches $G$: in this case this does not occur at a corner point. The right panel illustrates $\tau_1(\delta)$ in a similar manner, for $\delta =0.2,0.6$.}
 \label{fig:3dgausstau23}
\end{figure}

The general form of the gauge function for a meta-Gaussian distribution with standard exponential margins and correlation matrix $\Sigma$ with non-negative entries is 
\[
 g(\bm{x}) = (\bm{x}^{1/2})^\top \Sigma^{-1} \bm{x}^{1/2}.
\]
Figure~\ref{fig:3dgausstau23} displays the level set $g(\bm{x}) = 1$ when the Gaussian correlations in $\Sigma$ are $\rho_{12}=0.75, \rho_{13}=0.25, \rho_{23} =0.4$. The red dots on the level set are the points $(1,1,\gamma)/g(1,1,\gamma)$, $(1,\gamma,1)/g(1,\gamma,1)$ and $(\gamma,1,1)/g(\gamma,1,1)$ for $\gamma \in [0,1]$. The figure also provides an illustration of $\tau_{2,3}(\delta)$ for $\delta =0.2$ and $\delta=0.8$: in each case the light blue line from the origin is $\gamma \times (\delta,1,1)$, $\gamma \in [0,1]$, whilst the pink lines trace out the boundary $B_{\{2,3\},\delta}$ and $\tau_{2,3}(\delta)B_{\{2,3\},\delta}$. We see that when $\delta = 0.2$ (left panel), $\tau_{2,3}(0.2) = 1/g(0.2,1,1)$, i.e.,\ $\min_{\bm{y} \in B_{\{2,3\},0.2}}g(\bm{y}) = g(0.2,1,1)$. However, when $\delta = 0.8$, $\min_{\bm{y} \in B_{\{2,3\},0.8}}g(\bm{y}) = g(\gamma^\star,1,1)$, for $\gamma^\star \in [0,0.8]$, so $\tau_{2,3}(0.8) = 1/g(\gamma^\star,1,1)$. We note that the same value of $\tau_{2,3}(\delta)$ applies for any $\delta \geq \gamma^\star$: for this example, when $\delta \geq \gamma^\star \approx 0.51$, $\tau_{2,3}(\delta) = 0.7 = \eta_{2,3}$.

The reason that $\tau_{2,3}(\delta) = \eta_{2,3}$ for sufficiently large $\delta$ is because in this case $\arg\min_{x_1} g(x_1,1,1) = \gamma^\star$, meaning that the two-dimensional marginalization $g_{\{2,3\}}(1,1) = g(\gamma^\star,1,1)$, and we further have that $g_{\{2,3\}}(1,1) = \min_{\bm{y} \in B_{2,3}} g_{\{2,3\}}(\bm{y})$, so $\eta_{2,3} = 1/g_{2,3}(1,1)$. In Section~\ref{sec:vine} we will illustrate a gauge function for which $\arg\min_{x_3} g(1,1,x_3)>1$, and consequently $\tau_{1,2}(\delta) < \eta_{1,2}$ for all $\delta\leq 1$.

The right panel of Figure~\ref{fig:3dgausstau23} illustrates $\tau_{1}(\delta)$ for $\delta = 0.2$ and $\delta =0.6$. When $\delta =0.6$, the boundary $B_{1,\delta}$ already touches $G$, and so $\tau_{1}(0.6)=1$. In this example, $\tau_1(\delta) =1$ for any $\delta \geq 0.5625 = \rho_{12}^2$. As such, $\tau_{1}(0.2)<1$ as illustrated in the figure. We comment that if we had marginalized over $X_2$, and were looking at $\tau_1(\delta)$ for the variables $(X_1,X_3)$, then we would have $\tau_1(\delta) = 1$ for any $\delta \geq 0.0625 = \rho_{13}^2$. This provides an illustration of the dimensionality of the problem interacting with $\tau_{C}(\delta)$, and is again related to the point at which the minimum point defining the lower-dimensional gauge function occurs.

\subsubsection{Vine copula}
\label{sec:vine}
Three-dimensional vine copulas are specified by three bivariate copulas: two in the ``base layer'', giving the dependence between, e.g.,\ $X_1,X_2$ and $X_2,X_3$ and a further copula specifying the dependence between $X_1|X_2$ and $X_3|X_2$. Here we take the base copulas to be independence for $(X_1,X_2)$, and the inverted Clayton copula with parameter $\beta>0$ for $(X_2,X_3)$. The final copula is taken as inverted Clayton with parameter $\gamma>0$. The gauge function that arises in exponential margins is
\begin{align}
\label{eq:gvine}
g(\bm{x}) =  (1+\beta)\max(x_2,x_3)-\beta \min(x_2,x_3)&-\gamma x_1-(\gamma+1)(\beta+1)(\max(x_2,x_3)-x_2) \notag \\&+ (2\gamma+1)\max(x_1,(\beta+1)(\max(x_2,x_3)-x_2)).
\end{align}

\begin{figure}[]
 \includegraphics[width=0.9\textwidth]{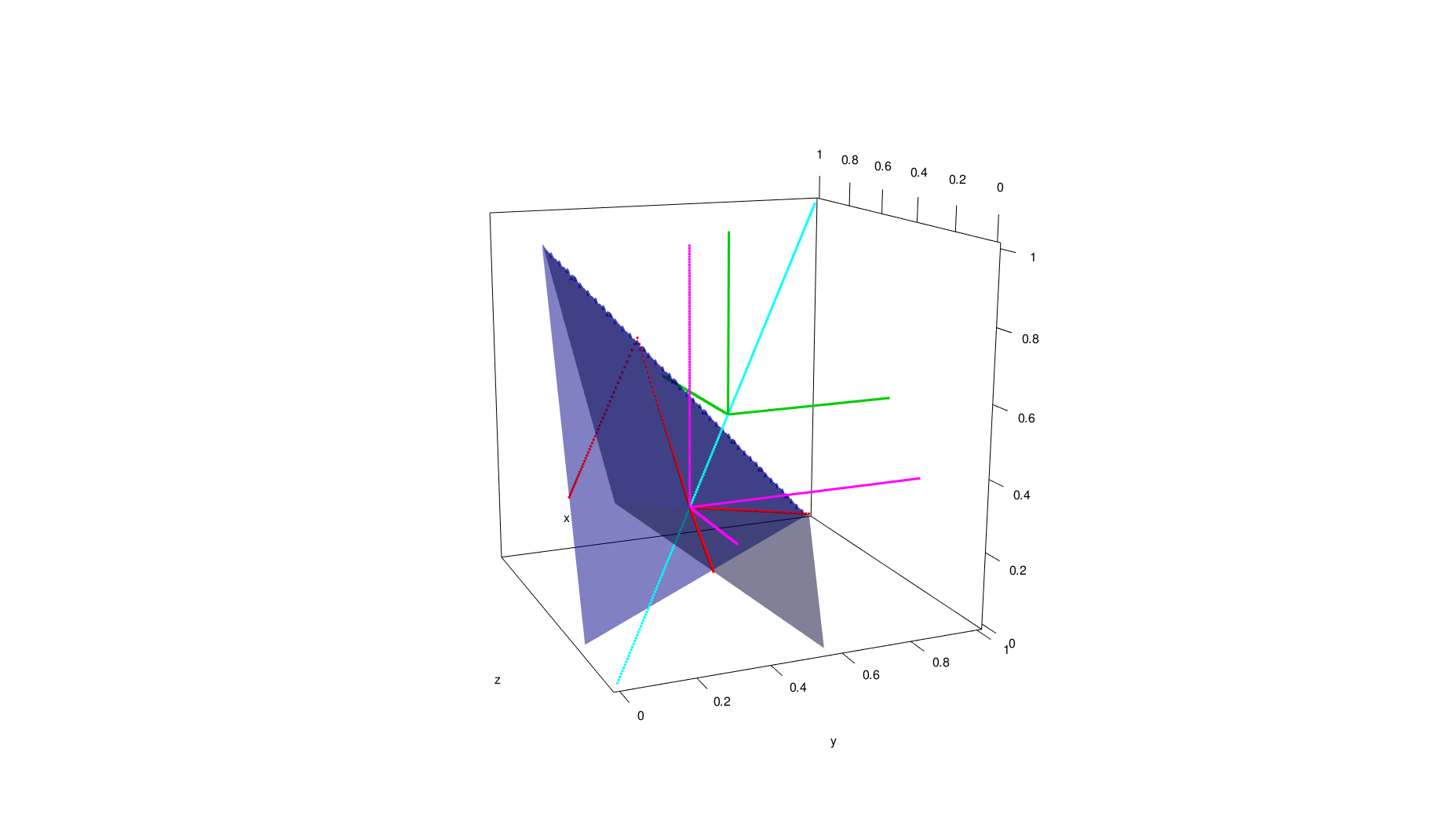}
 \caption{Level set $g(\bm{x})=1$ with $g$ as in~\eqref{eq:gvine}. The figure illustrates $\tau_{1,2}(1)$ and $\eta_{1,2,3} = \eta_{1,2} = 1/2$.}
 \label{fig:3dvinetau12}
\end{figure}

Figure~\ref{fig:3dvinetau12} displays the level set $g(\bm{x})=1$. In this figure we also give an illustration of a case where $\tau_C(1)<\eta_C$: in particular, for this example $\tau_{1,2}(1)<\eta_{1,2,3}=\eta_{1,2}=1/2$. The purple lines represent the boundary of the region $\tau_{1,2}(1)R_{\{1,2\},1} = \tau_{1,2}(1,\infty]^2 \times [0,1]$, while the green lines represent the boundary of the region $\eta_{1,2,3}(1,\infty]^3$. Theorem~1 of \citet{Simpsonetal18} tells us that $\eta_{1,2} = \max(\tau_{1,2}(1),\tau_{1,2,3})$, where $\tau_{1,2,3}=\eta_{1,2,3}$. Therefore $\tau_{1,2}(1)<\eta_{1,2}$ guarantees that $\eta_{1,2}=\eta_{1,2,3}$.

We also illustrate Proposition~\ref{prop:marg}, minimizing~\eqref{eq:gvine} over $x_3$. If $x_2>x_3$ then the minimum over $x_3$ occurs by setting $x_3=x_2$ and is equal to $x_2+(1+\gamma)x_1$. If $x_2< x_3$ then owing to the final term we need to consider the cases $x_3 \lessgtr x_1/(1+\beta)+x_2$. In both cases, the mimimum is attained at $x_3 = x_1/(1+\beta)+x_2$, and is equal to $x_1+x_2 < (1+\gamma)x_1+x_2$. As such, $\min_{x_3} g(x_1,x_2,x_3) = x_1+x_2$. This result is as expected since the bivariate margins of vine copulas that are directly specified in the base layer are equal to the specified copula: in this case, independence.

\subsection{Mixing independent vectors}
\label{sec:mixing}

Here we exploit the results from previous sections to consider what happens when independent exponential random vectors are additively mixed such that the resulting vector still has exponential type tails. We consider as a case study the spatial model of \citet{HuserWadsworth18}, which following a reparameterization can be expressed
\begin{align}
    \{X_{\tilde{E}}(s) = \gamma S_E + V_E(s): s \in \mathcal{S}\subset \mathbb{R}^2\}, \qquad \gamma \in (0,\infty), \label{eq:HWspat}
\end{align}
where $S_E \sim \mbox{Exp}(1)$ is independent of the spatial process $V_{E}$, which also possesses unit exponential margins and is asymptotically independent at all spatial lags $s_1-s_2 \neq 0$. The process $V_E$ is assumed to possess hidden regular variation, with residual tail dependence coefficient satisfying $\eta^{V}(s_1,s_2)<1$ for all $s_1 \neq s_2$. The resulting process $X_{\tilde{E}}$ is asymptotically independent for $\gamma \in (0,1]$ and asymptotically dependent for $\gamma>1$; see also \citet{Engelkeetal18} for related results.

When $\gamma <1$, $\pbb(X_{\tilde{E}}(s) > x) \sim e^{-x}/(1-\gamma)$. In this case, \citet{HuserWadsworth18} show that the residual tail dependence coefficient for the process $X_{\tilde{E}}$ is given by
\begin{align}
    \eta^{X} = \begin{cases}
    \eta^{V}, & \gamma < \eta^V\\
    \gamma, & \eta^V \leq \gamma \leq 1. 
    \end{cases} \label{eq:etaHW}
\end{align}
That is, the strength of the extremal dependence as measured by the residual tail dependence coefficient $\eta^X$ is increasing in $\gamma$ for $\gamma \geq \eta^V$. In contrast, \citet{WadsworthTawn19} showed that under mild conditions, the process~\eqref{eq:HWspat} has the same conditional extremes normalization as the process $V_E(s)$, with identical limit distribution when the scale normalizations $b_{s-s_0}(t) \to \infty$ as $t \to \infty$. Here, the subscript $s-s_0$ alludes to the fact that the conditioning event in~\eqref{eq:ce} is $\{V_{E}(s_0)>t\}$ and we study the normalization at some other arbitrary location $s \in \mathcal{S}$. In combination, we see that the results of \citet{HuserWadsworth18} and \citet{WadsworthTawn19} suggest that the addition of the variable $\gamma S_E$ to $V_E$ affects the extremal dependence of $X_{\tilde{E}}$ differently for different extreme value representations. We elucidate these results further in the context of the limit sets and their gauge functions. \iftoggle{arxiv}{}{A summary is provided here, with full derivations in the arXiv version.}

Let us suppose that $\bm{S}_E \in \mathbb{R}^{d}_+$ has unit exponential margins, density $f_{S_E}$, and gauge function $g_S$, and is independent of $\bm{V}_E \in \mathbb{R}^{d}_+$, which has unit exponential margins, density $f_{V_E}$, and gauge function $g_V$. Let $\bm{Z}_E = (\bm{S}_E,\bm{V}_E) \in \mathbb{R}^{2d}_+$ be the concatenation of these vectors\iftoggle{arxiv}{. Then, since its density $f_{Z_{E}}(\bm{z}) = f_{S_E}(z_1,\ldots,z_{d}) f_{V_E}(z_{d+1},\ldots,z_{2d})$, it is clear that $\bm{Z}_E$ has unit exponential margins and gauge function $g_Z(\bm{z}) = g_{S}(z_1,\ldots,z_{d}) + g_{V}(z_{d+1},\ldots,z_{2d})$.}{; this has exponential margins and gauge function $g_Z(\bm{z}) = g_{S}(z_1,\ldots,z_{d}) + g_{V}(z_{d+1},\ldots,z_{2d})$.}

Now consider the linear transformation of $\bm{Z}_E$ to 
\begin{align*}
 A\bm{Z}_E =  (\gamma Z_{E,1}+ Z_{E,d+1},\ldots, \gamma Z_{E,d}+ Z_{E,2d}, Z_{E,1},\ldots, Z_{E,d}) = (\gamma \bm{S}_E+\bm{V}_{E}, \bm{S}_E) = (\bm{X}_{\tilde{E}},\bm{S}_E),
\end{align*}
where $A \in \mathbb{R}^{2d \times 2d}$ is the matrix describing this transformation\iftoggle{arxiv}{: the first $d$ rows have $\gamma$ and $1$ in the $j$th and $j+d$th positions, respectively, for $j=1,\ldots,d$, while the second $d$ rows have $1$ in the $j$th position for $j=1,\ldots, d$. All other entries are zero. The matrix $A^{-1}$ has the same configuration but with $-\gamma$ in place of $\gamma$.}{.}
By Lemma~4.1 of \citet{Nolde14}, the normalized sample cloud $\{A\bm{Z}_{E,i}/\log n : i=1,\ldots, n\}$ converges onto the set $AG$, where $G=\{\bm{z}\in\mathbb{R}^{2d}_+:g_Z(\bm{z}) \leq 1\}$, so $AG = \{\bm{z}\in\mathbb{R}^{2d}_+:A^{-1}\bm{z} \in G\} = \{\bm{z}\in\mathbb{R}^{2d}_+: g_Z(A^{-1}\bm{z}) \leq 1\}$. Consequently, the gauge function of $A\bm{Z}_E$ is $g_Z(A^{-1}\bm{z})$, i.e., $g_{Z}(\bm{x},\bm{s})=g_{S}(\bm{s}) + g_{V}(\bm{x}-\gamma\bm{s})$, for $\bm{x}>\gamma\bm{s}$.

Next we apply Proposition~\ref{prop:marg} to the vector $A\bm{Z}$, marginalizing over the last $d$ coordinates, which are equal to $\bm{S}_{E}$. This leaves us with the gauge function of $X_{\tilde{E}}$, denoted $g_X$, and given by
\begin{align*}
    g_X(\bm{x}) = \min_{\bm{s} \in [\bm{0},\bm{x}/\gamma]} g_S(s_1,\ldots,s_d) + g_{V}(x_1-\gamma s_1,\ldots, x_d - \gamma s_d).
\end{align*}

To illustrate the results of \citet{HuserWadsworth18} and \citet{WadsworthTawn19} concerning model~\eqref{eq:HWspat}, we need to take $\bm{S}_{E} = S_{E}\bm{1}$, i.e., perfect dependence. Although such a vector does not have a $d$-dimensional Lebesgue density, convergence of the sample cloud based on the univariate random variable $S_E$ onto the unit interval $[0,1]$ implies that the limit set is $G_S = \{\bm{x} \in \mathbb{R}^d_+: x_1=x_2=\cdots=x_d=x, x \leq 1\}$. Such a set can be described by the gauge function 
\[
g_S(\bm{s}) = \begin{cases}
\infty, & s_i \neq s_j~~\mbox{for any}~~i,j\\
s, & s_1=\cdots=s_d = s.
\end{cases}
\]
As such, in this case, $ g_X(\bm{x}) = \min_{s \in[0,\min(\bm{x})/\gamma]} \{s + g_{V}(\bm{x}-\gamma s)\}.$

\paragraph{Residual tail dependence} To find the residual tail dependence coefficient $\eta^X$, we require
\begin{align*}
    \min_{\bm{x}:\min(\bm{x})=1} g_{X}(\bm{x}) &=   \min_{\bm{x}:\min(\bm{x})=1}\min_{s \in[0,\min(\bm{x})/\gamma]} \{s + g_{V}(\bm{x}-\gamma s)\}\\
    &=\min_{s \in[0,1/\gamma]}\min_{\bm{x}:\min(\bm{x})=1} \{s + g_{V}(\bm{x}-\gamma s)\}.
\end{align*}
For fixed $s$, consider $\min_{\bm{x}:\min(\bm{x})=1} g_{V}(\bm{x}-\gamma s) = \min_{\bm{z}:\min(\bm{z})=1-\gamma s} g_{V}(\bm{z}) = g_{V}(\bm{y}^\star\times(1-\gamma s))$, where $\bm{y}^\star = \arg\min_{\bm{y}:\min(\bm{y})=1} g_V(\bm{y})$. As such
\begin{align*}
    \min_{\bm{x}:\min(\bm{x})=1} g_{X}(\bm{x}) &= \min_{s \in[0,1/\gamma]}\{s + g_V(\bm{y}^\star)(1-\gamma s)\} = \begin{cases} g_V(\bm{y}^\star),& \gamma < 1/g_V(\bm{y}^\star),\\ 1/\gamma, & \gamma \geq 1/g_V(\bm{y}^\star). \end{cases}
\end{align*}
Recalling that $\eta^X = [\min_{\bm{x}:\min(\bm{x})=1} g_{X}(\bm{x})]^{-1}$ and $\eta^V=1/g_{V}(\bm{y}^\star)$, this yields~\eqref{eq:etaHW}.

\paragraph{Conditional extremes} For the conditional extremes normalization, we now let $g_V$ and $g_X$ denote two-dimensional gauge functions. Suppose that $\alpha_V, \beta_V$ are such that $g_V(\alpha_V,1) = 1$ and $g_V(\alpha_V + u,1) - 1 \in \RV_{1/(1-\beta_V)}^0$. We have 
\begin{align}
   1= g_X(\alpha_X,1) = \min_{s} \{s + g_{V}(\alpha_X-\gamma s,1-\gamma s)\}. \label{eq:ceg}
\end{align}
Suppose that the right hand side of~\eqref{eq:ceg} is minimized at $s^\star \geq 0$, i.e., $g_X(\alpha_X,1) = s^\star + g_{V}(\alpha_X-\gamma s^\star,1-\gamma s^\star)$. Because $\alpha_X \leq 1$ and $g_V(v_1,v_2)\geq \max(v_1,v_2)$, this yields 1=$g_X(\alpha_X,1) \geq 1+ (1-\gamma) s^\star$, therefore we must have $s^\star=0$ for $\gamma \in (0,1)$. Consequently, $\alpha_X = \alpha_V = \alpha$. 

\iftoggle{arxiv}{For the scale normalization, let $s_u^\star = \arg\min_{s \geq 0} \{s + g_{V}(\alpha_X-\gamma s + u,1-\gamma s)\}$. Then 
\begin{align}
   0 \leq  g_{X}(\alpha + u, 1) -1 = s_u^\star + g_{V}(\alpha-\gamma s_u^\star + u,1-\gamma s_u^\star) - 1 \leq g_{V}(\alpha + u, 1) -1, \label{eq:ineqst1}
\end{align}
and because $g_V(x,y) \geq \max(x,y)$, if $\alpha<1$ then for sufficiently small $u$
\begin{align}
   (1-\gamma)s_u^\star \leq s_u^\star + g_{V}(\alpha-\gamma s_u^\star + u,1-\gamma s_u^\star) - 1 \label{eq:ineqst2}.
\end{align}
Combining inequalities~\eqref{eq:ineqst1} and~\eqref{eq:ineqst2} shows that for $\gamma \in (0,1)$
\begin{align}
0 \leq (1-\gamma)s_u^\star \leq  g_{X}(\alpha + u, 1) -1 \leq g_{V}(\alpha + u, 1) - 1 \to 0, \qquad u \to 0, \label{eq:ineqst3}
\end{align}
meaning in particular that $s_u^\star \to 0$ as $u \to 0$. To examine the minimizing sequence $s_u^\star$ in further detail, we consider the derivative of $s + g_{V}(\alpha-\gamma s + u,1-\gamma s)$, assuming that $g_V$ has piecewise continuous partial derivatives possessing finite left and right limits, denoted by $g_{V,1}, g_{V,2}$, such that
\begin{align}
\frac{\partial}{\partial s} \left[s + g_{V}(\alpha-\gamma s + u,1-\gamma s) \right]& = 1-\gamma\left[g_{V,1}(\alpha+u -\gamma s,1-\gamma s)+ g_{V,2}(\alpha+u -\gamma s,1-\gamma s)\right] \notag \\
& = 
1-\gamma\left[g_{V}(\alpha+u -\gamma s,1-\gamma s) + \{1-(\alpha+u - \gamma s)\} g_{V,1}(\alpha+u -\gamma s,1-\gamma s)\right. \notag \\ & \qquad\qquad \left. -\gamma s g_{V,2}(\alpha+u -\gamma s,1-\gamma s)\right], \label{eq:gXderiv}
\end{align}
using Euler's homogeneous function theorem on the second line. Consider~\eqref{eq:gXderiv} evaluated at $s_u^\star$ and take the limit as $u \to 0$, yielding
\begin{align}
1-\gamma[1+(1-\alpha)g_{V,1}(\alpha_+,1)]. \label{eq:gXderivlim}
\end{align}
For $g_V$ differentiable at $(\alpha,1)$, $g_{V,1}(\alpha_+,1) = g_{V,1}(\alpha_-,1) = 0$, so the limit of~\eqref{eq:gXderiv} is $1-\gamma>0$, and hence there exists $\epsilon>0$ such that \eqref{eq:gXderiv} is positive for all $u< \epsilon$, giving $s_u^\star = 0$ for all $u<\epsilon$. Consequently, $g_{X}(\alpha+u,1)-1 \sim g_{V}(\alpha+u,1)-1$, $u \to 0$.

For $g_{V}$ not differentiable at $(\alpha,1)$,~\eqref{eq:gXderivlim} is positive when $g_{V,1}(\alpha_+,1)<(1/\gamma-1)/(1-\alpha)$, and again in this case $g_{X}(\alpha+u,1)-1 \sim g_{V}(\alpha+u,1)-1$, $u \to 0$. When $g_{V,1}(\alpha_+,1)>(1/\gamma-1)/(1-\alpha)$, then the minimizing sequence $s_u^\star$ should be as large as possible, i.e., equal to its upper bound of $[g_{X}(\alpha+u,1)-1]/(1-\gamma) = [s_u^\star + g_{V}(\alpha-\gamma s_u^\star + u,1-\gamma s_u^\star) - 1]/(1-\gamma)$ from inequality~\eqref{eq:ineqst3}. Further asymptotic detail on this bound is obtained through a Taylor expansion:
\begin{align}
 g_{V}(\alpha-\gamma s_u^\star + u,1-\gamma s_u^\star) = g_V(\alpha,1) + (u-\gamma s_u^\star)g_{V,1}(\alpha_+,1) - \gamma s_u^\star g_{V,2}(\alpha,1_-) + O(\max(u^2,s_u^{\star 2},u s_u^\star)), \label{eq:gXTaylor}
\end{align}
giving
\begin{align*}
s_u^\star \leq \frac{u g_{V,1}(\alpha_+,1)}{\gamma[g_{V,1}(\alpha_+,1)+g_{V,2}(\alpha,1_-) - 1]}[1+o(1)] = \frac{u}{\gamma(1-\alpha)} + o(u).
\end{align*}
Taking $s_u^\star$ at this upper bound and using the expansion~\eqref{eq:gXTaylor}, 
\begin{align*}
g_{X}(\alpha+u,1)-1 = s_u^\star + g_{V}(\alpha-\gamma s_u^\star + u,1-\gamma s_u^\star) -1 = \frac{1-\gamma}{\gamma(1-\alpha)}u + o(u),
\end{align*}
such that $g_{X}(\alpha+u,1)-1 \in \RV^0_{1}$. Since we also have $g_{V}(\alpha+u,1)-1 \sim u g_{V,1}(\alpha_+,1)\in \RV^0_{1}$, the regular variation indices are identical with $\beta_X=\beta_V=0$. This represents a case where the scale normalizations $b_{s-s_0}(t)$ in the conditional extremes representation do not diverge to infinity, meaning a potential difference in the limit distribution.}{Calculations for the scale normalization are more involved and can be found in the arXiv version. We find that for $g_V$ differentiable at $(\alpha,1)$, $g_{X}(\alpha+u,1)-1 \sim g_{V}(\alpha+u,1)-1$, $u \to 0$, whereas in the non-differentiable case we do not necessarily have this link but can deduce that the regular variation indices are $\beta_X=\beta_V = 0$.}

\begin{figure}
    \centering
    \includegraphics[width=0.3\textwidth]{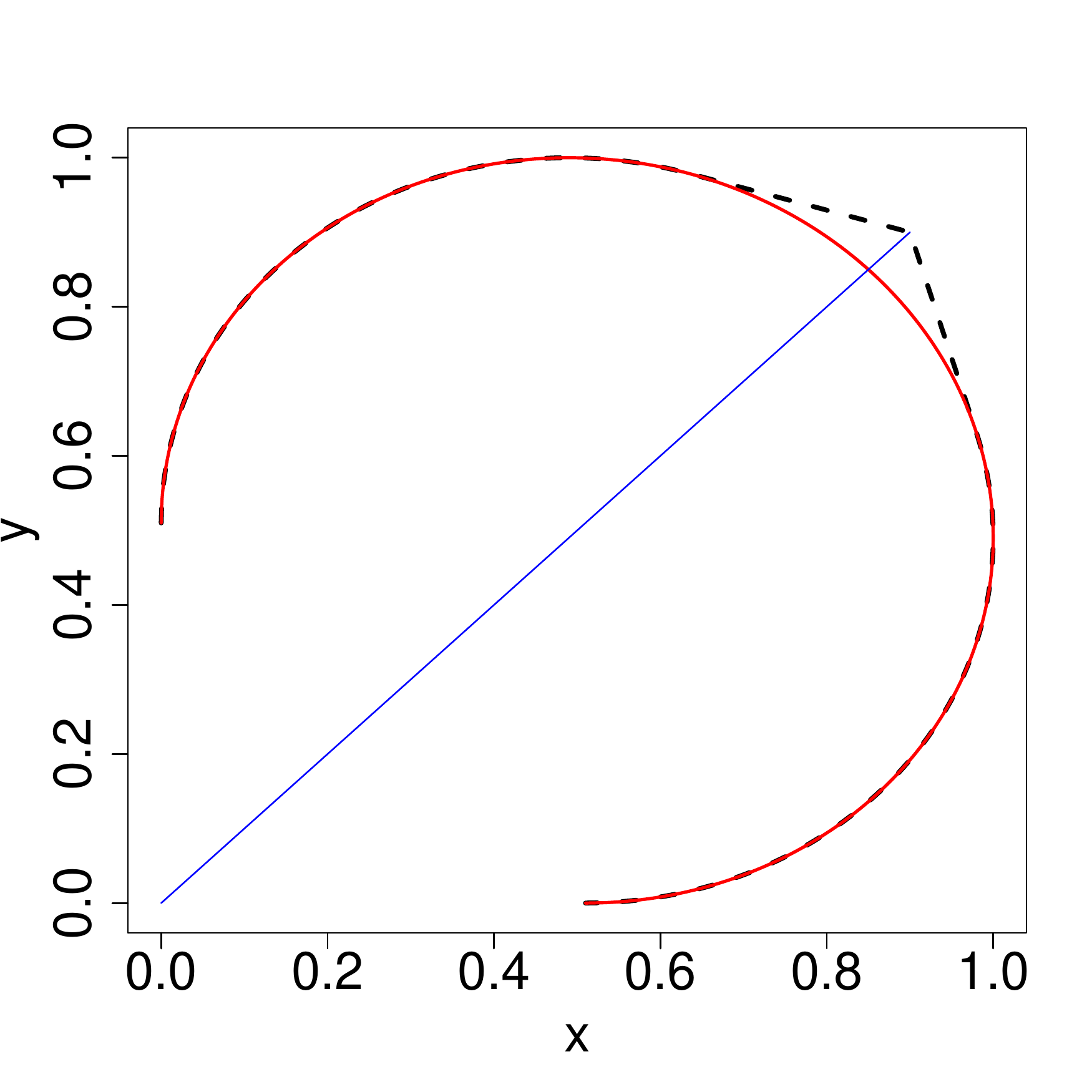}
    \includegraphics[width=0.3\textwidth]{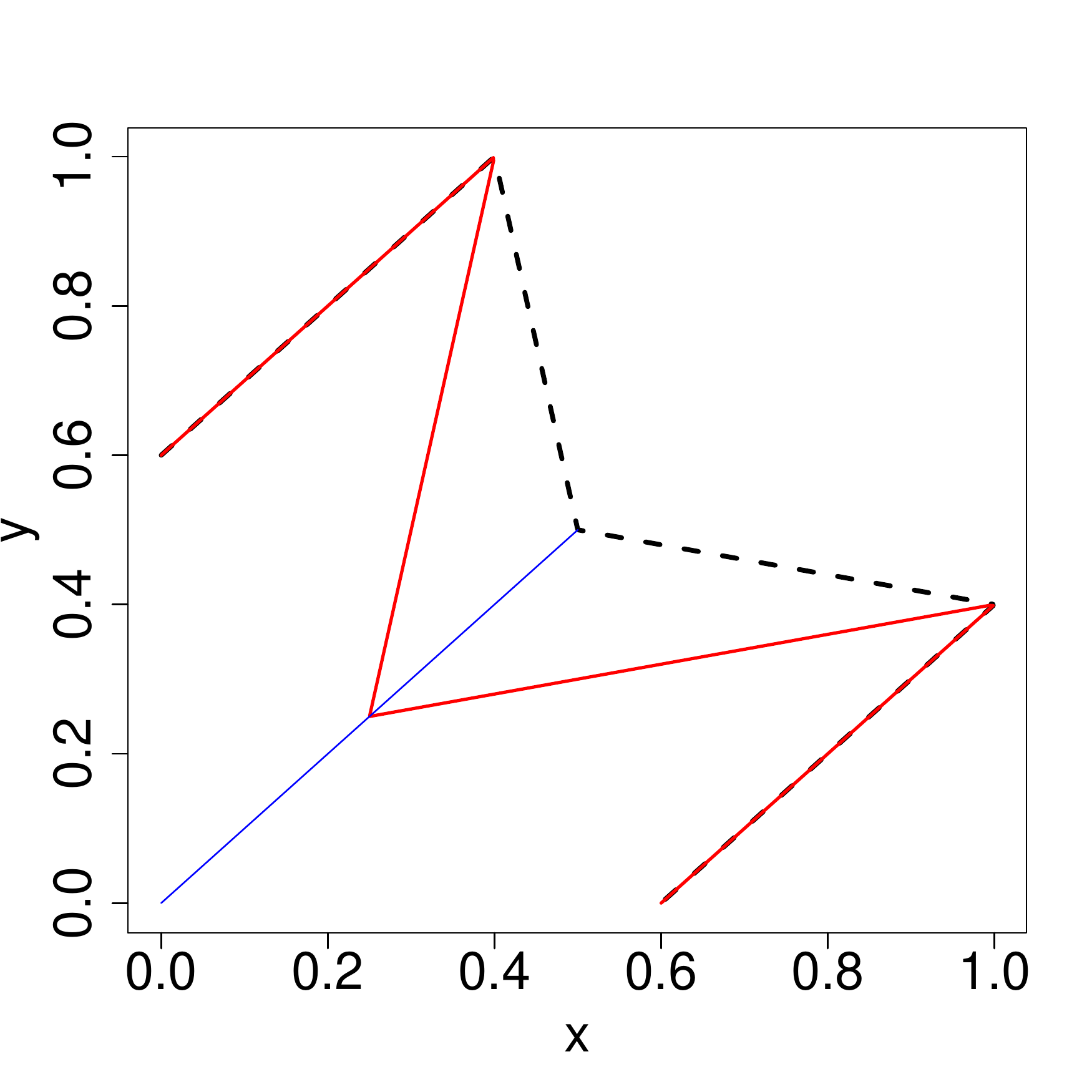}
    \includegraphics[width=0.3\textwidth]{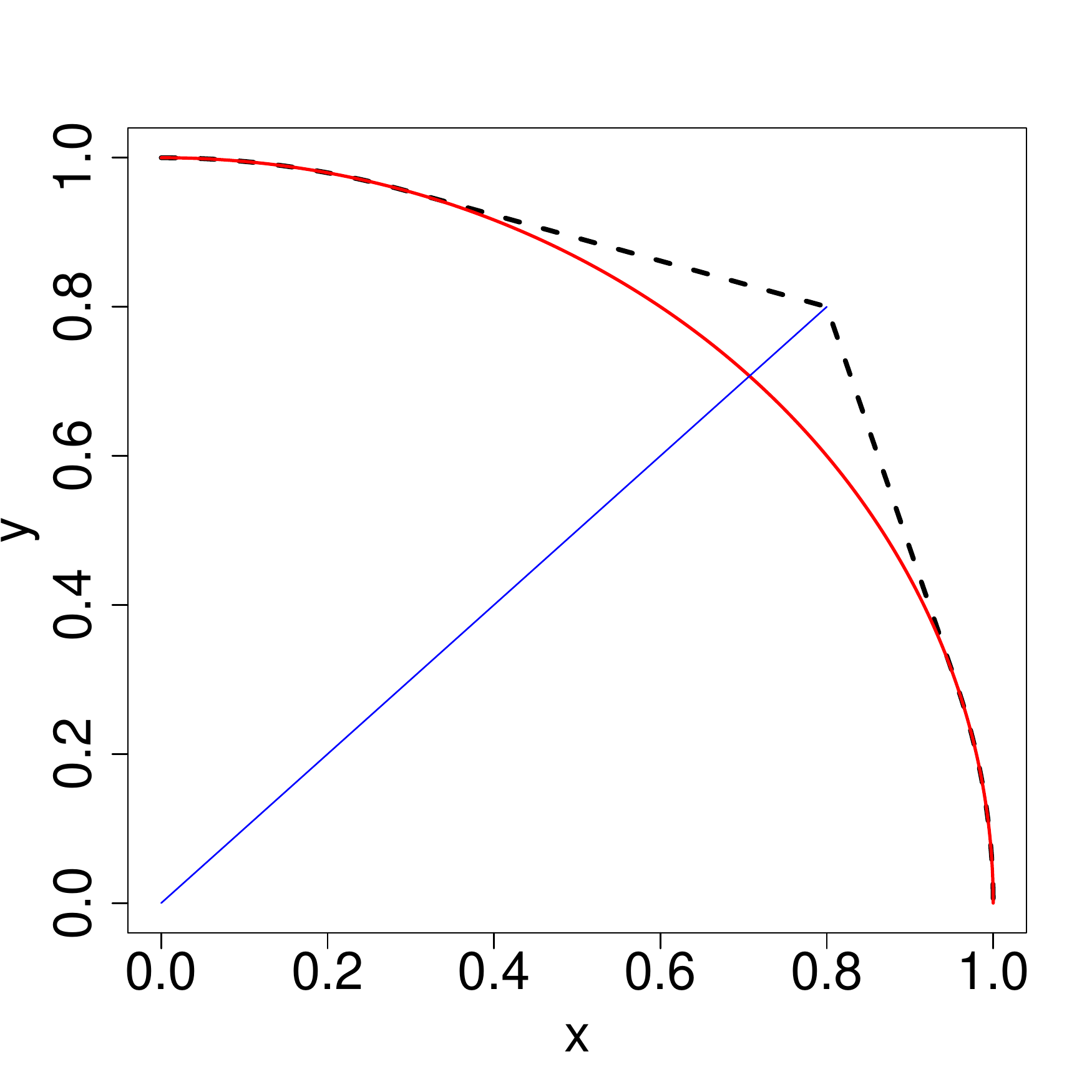}
    \caption{Solid red lines depict the level sets $g_V(\bm{x})=1$, where $g_V$ is of the form (ii), (iii) and (iv) (L--R) from Figure~\ref{fig:2dgauge}. Dashed black lines depict the level sets $g_X(\bm{x}) =1$. In each picture the blue solid line segment is from $(0,0)$ to $(\gamma,\gamma)$, and denotes the limit set of the fully dependent random vector $\gamma S_E$. From left to right, $\gamma=0.9,0.5,0.8$.}
    \label{fig:HW}
\end{figure}

\iftoggle{arxiv}{Figure~\ref{fig:HW} displays examples of gauge functions $g_V$ and $g_X$. We observe from this figure how, when $\gamma$ is sufficiently large, the shape of $g_V$ is modified to produce $g_X$. The modification is focussed around the diagonal, and explains visually why the residual tail dependence coefficient changes while the conditional extremes normalization does not. The left and right panels illustrate differentiable cases, with $s_u^\star \equiv 0$ for sufficiently small $u$; the centre panel depicts an example with $s_u^\star$ linear in $u$ as $u\to 0$.}{Figure~\ref{fig:HW} displays examples of gauge functions $g_V$ and $g_X$. We observe from this figure how, when $\gamma$ is sufficiently large, the shape of $g_V$ is modified to produce $g_X$. The modification is focussed around the diagonal, and explains visually why the residual tail dependence coefficient changes while the conditional extremes normalization does not. The left and right panels illustrate differentiable cases, and the centre panel non-differentiable.}

\section{Discussion}
\label{sec:Discussion}
In this work we have demonstrated how several concepts of extremal dependence can be unified through the shape of the limit set $G$ of the scaled sample cloud $N_n =\{\bX_1/r_n,\ldots,\bX_n/r_n\}$ arising for distributions with light-tailed margins. For concreteness our focus has been on exponential margins, but other choices can be useful. In the case of negative dependence between extremes --- such that large values of one variable are most likely to occur with small values of another --- the double exponential-tailed Laplace margins can be more enlightening. As an example, for the bivariate Gaussian copula with $\rho<0$ we observed that the limit set $G$ is described by a discontinuous gauge function $g$ that cannot be established through the simple mechanism of Proposition~\ref{prop:logdens}. In \citet{Nolde14}, the gauge function for this distribution in Laplace margins is calculated as
\begin{align*}
g(x,y) = \begin{cases}
(|x|+|y|-2\rho|xy|^{1/2})/(1-\rho^2), & x,y \geq 0~~\mbox{or}~~x,y \leq 0, \\
(|x|+|y|+2\rho|xy|^{1/2})/(1-\rho^2), & x\geq 0, y \leq 0~~\mbox{or}~~x\leq 0,y \geq 0.
\end{cases}
\end{align*}
When $\rho<0$, this yields $g(1,-\rho^2) = 1$ , and $g(1,-\rho^2 + u) \in \RV_{2}^0$, so that extending Proposition~\ref{prop:CEgauge}, we would find that the conditional extremes normalizations are $a^j(t) \sim -\rho^2 t$ and $b^j(t) \in \RV^0_{1/2}$, as given in \citet{Keefetal13}.

The study of extremal dependence features through the limit set $G$ is enlightening both for asymptotically dependent and asymptotically independent random vectors, particularly as it can be revealing for mixture structures where mass is placed on a variety of cones $\bbE_C$ as defined in~\eqref{eq:EC}. However, many traditional measures of dependence within the asymptotically dependent framework, which are typically functions of the exponent function $V$ given in equation~\eqref{eq:V}, or spectral measure $H$, are not revealed by limit set~$G$. For example, it was noted in the example of Section~\ref{sec:eglog} that the limit set described by the  gauge function $g(x,y) = \theta^{-1}\max(x,y) + (1-\theta^{-1})\min(x,y)$ can arise for several different spectral measures, although clearly the parameter $\theta$ demonstrates some link between strength of dependence and shape of $G$.

 Nonetheless, multivariate regular variation and associated limiting measures have been well-studied in extreme value theory, but representations that allow greater discrimination between asymptotically independent or mixture structures much less so. The limit set elucidates many of these alternative dependence concepts and provides meaningful connections between them.  We have not directly considered connections between the various dependence measures without reference to $G$, and we note that the limit set might not always exist. We leave such study to future work.

\subsection*{Acknowledgements}
The authors are grateful to the editor and two reviewers for constructive feedback and valuable comments that helped improve the paper. NN acknowledges financial support of the Natural Sciences and Research Council of Canada. JLW gratefully acknowledges funding from EPSRC grant EP/P002838/1.

\bibliographystyle{apalike}

\bibliography{GaugeBib}

\end{document}